\documentclass[11pt]{amsart}
\usepackage[utf8]{inputenc}
\usepackage{comment}
\usepackage[style=alphabetic,firstinits,backend=bibtex, uniquename=init,
url=true,hyperref]{biblatex}
\DeclareFieldFormat{labelalpha}{\thefield{entrykey}}
\DeclareFieldFormat{extraalpha}{}
\DeclareFieldFormat{postnote}{#1}
\DeclareFieldFormat{multipostnote}{#1}
\usepackage{tikz}
\usetikzlibrary{decorations.text,calc,arrows.meta}
\usepackage{booktabs}
\usepackage{multirow}
\usepackage{array}
\usepackage{varwidth}
\usepackage{filecontents}
\usepackage[titletoc]{appendix}
\usepackage{mathtools}
\usepackage{amsthm}
\usepackage{amsmath}
\usepackage{amsfonts}
\usepackage{amssymb}
\usepackage{mathrsfs}
\usepackage{xypic}
\usepackage{color}
\usepackage{xcolor}
\usepackage{titlesec}
\usepackage{titletoc}
\usepackage[colorlinks]{hyperref}
\definecolor{cclr}{rgb}{25,25,112}
\hypersetup{
citecolor=[rgb]{0.15, 0.15, 0.68},
linkcolor=gray,   
urlcolor=magenta}

\bibliography{eg-l-par}

\newtheorem{lem}{Lemma}
\newtheorem{sublem}{Sublemma}
\newtheorem{cor}{Corollary}
\newtheorem{prop}{Proposition}
\newtheorem{exam}{Example}
\newtheorem{thm}{Theorem}

\newtheorem{remark}{Remark}
\newtheorem{dfn}{Definition}
\numberwithin{thm}{section}
\numberwithin{lem}{section}
\numberwithin{cor}{section}
\numberwithin{prop}{section}
\numberwithin{exam}{section}
\numberwithin{dfn}{section}
\def\presuper#1#2%
  {\mathop{}%
   \mathopen{\vphantom{#2}}^{#1}%
   \kern-\scriptspace%
   #2}

\newcommand{\sep}{{\operatorname{sep}}}
\newcommand{\Inf}{{\operatorname{inf}}}

\newcommand{\wh}[1]{\widehat{#1}}
\newcommand{\wt}[1]{\widetilde{#1}}
\newcommand{\Ad}{\operatorname{Ad}}

\newcommand{\Lie}{\operatorname{Lie}}
\newcommand{\Sym}{\operatorname{Sym}}

\newcommand{\End}{\operatorname{End}}

\newcommand{\std}{\operatorname{std}}

\newcommand{\Gal}{\operatorname{Gal}}

\newcommand{\red}{{\operatorname{red}}}

\newcommand{\frep}{^f\operatorname{Rep}}

\newcommand{\frepg}{^f\operatorname{Rep}_G}
\newcommand{\gtor}{\cG\operatorname{Tor}}

\newcommand{\vect}{\operatorname{Vect}}

\newcommand{\Aut}{\operatorname{Aut}}
\newcommand{\id}{\operatorname{id}}
\newcommand{\Res}{\operatorname{Res}}
\newcommand{\univ}{{\operatorname{univ}}}

\newcommand{\Herr}{{\operatorname{Herr}}}
\newcommand{\GR}{\operatorname{Gr}}

\newcommand{\GL}{\operatorname{GL}}

\newcommand{\Img}{\operatorname{Im}}

\newcommand{\Pt}{\operatorname{Pt}}
\newcommand{\Vect}{\operatorname{Vect}}
\newcommand{\QCoh}{\operatorname{QCoh}}
\newcommand{\Coh}{\operatorname{Coh}}
\newcommand{\Hom}{\operatorname{Hom}}

\newcommand{\iHom}{\underline{\operatorname{Hom}}}

\newcommand{\spf}{\operatorname{Spf}}
\newcommand{\spec}{\operatorname{Spec}}

\newcommand{\invlim}[1]{{\varprojlim_{#1}}}
\newcommand{\dirlim}[1]{{\varinjlim_{#1}}}
\titleformat{\section}[runin]{\normalfont\bfseries}{\thesection.}{3pt}{}
\titleformat{\subsection}[runin]{\normalfont\bfseries}{\thesubsection.}{3pt}{}
\titleformat{\subsubsection}[runin]{\normalfont\bfseries}{\thesubsubsection.}{3pt}{}
\titleformat{\paragraph}[runin]{\normalfont\bfseries}{}{3pt}{}
\renewcommand{\thesection}{\arabic{section}}
\newcommand{\mychapter}{\section}
\newcommand{\mysection}{\subsection}
\newcommand{\mysubsection}{\subsubsection}
\titleformat{\section}{\normalfont\large\bfseries}{\thesection.~~}{1em}{}
 
\newcommand{\fc}{\mathfrak{c}}

\newcommand{\cO}{\mathcal{O}}
\newcommand{\fS}{\mathfrak{S}}
\newcommand{\cK}{\mathcal{K}}
\newcommand{\fm}{\mathfrak{m}}
\newcommand{\fp}{\mathfrak{p}}
\newcommand{\fq}{\mathfrak{q}}
\newcommand{\cC}{\mathcal{C}}
\newcommand{\fC}{\mathfrak{C}}

\newcommand{\fT}{\mathfrak{T}}

\newcommand{\bC}{\mathbb{C}}
\newcommand{\cX}{\mathcal{X}}
\newcommand{\cY}{\mathcal{Y}}
\newcommand{\cZ}{\mathcal{Z}}
\newcommand{\cV}{\mathcal{V}}

\newcommand{\calR}{\mathcal{R}}
\newcommand{\clR}{\mathcal{R}}

\newcommand{\cF}{\mathcal{F}}

\newcommand{\cG}{\mathcal{G}}
\newcommand{\clH}{\mathcal{H}}
\newcommand{\G}{\mathbb{G}}
\newcommand{\calD}{\mathcal{D}}
\newcommand{\cE}{{\mathcal{E}}}
\newcommand{\et}{\text{\'et}}
\newcommand{\ef}{\text{ef}}
\newcommand{\bI}{\mathbb{I}}
\newcommand{\bFp}{\bar{\F}_p}
\newcommand{\bM}{\breve{M}}
\newcommand{\bP}{\breve{P}}
\newcommand{\bB}{\breve{B}}
\newcommand{\bT}{\breve{T}}
\newcommand{\bG}{\breve{G}}

\newcommand{\E}{\mathbf{E}}
\newcommand{\A}{\mathbf{A}}
\newcommand{\Z}{\mathbf{Z}}
\newcommand{\Q}{\mathbf{Q}}
\newcommand{\C}{\mathbf{C}}
\newcommand{\F}{\mathbf{F}}

\newcommand{\Qp}{\Q_p}
\newcommand{\lsup}[2]{{^{#1}\!#2}}
\newcommand{\lG}{{\lsup LG}}

\newcommand{\Kummer}{\operatorname{Kummer}}

\dottedcontents{section}[2.5em]{}{2.9em}{1pc}

\begin{document}
\author{Lin, Zhongyipan}
\title{Moduli stacks of generalized $\varphi$-modules}

\begin{abstract}
Let $\varphi\in \End(\Z_q/p^a(\!(u)\!))$
be a ring endomorphism such that
$\varphi(u)\equiv u^{p}$ mod $p$.
We show that,
as long as a {\it height theory} is definable,
the moduli stack
of \'etale $\varphi$-modules
(with $G$-structure)
of height {\it bounded} by a constant
is an algebraic stack of finite presentation
over $\Z/p^a$.
Most notably, we allow
$\varphi(u)\notin \Z_q/p^a[\![u]\!]$.

Let $F/\Q_p$ be an arbitrary $p$-adic field
and let $G$ be an arbitrary reductive group over $F$
with Langlands dual group $\lsup LG$.
It follows that the Emerton-Gee stack
$\cX_{F,\lsup LG}\cong\cX_{\Q_p,\lsup L\Res_{F/\Q_p}G}
\to\cX_{\Q_p, \GL_d}$
is relatively representable by {\it algebraic stacks}
of finite presentation over $\spf \Z_p$
for any embedding $\lsup L\Res_{F/\Q_p}G
\to\GL_d$,
which improves the result of \cite{Min25}
which says the morphism is representable by
locally Noetherian
{\it formal} algebraic stacks;
we also establish
a weaker representability result 
for parabolic subgroups $\lsup LP\subset\lsup LG$.
Meanwhile, the algebraicity
of moduli of $\varphi$-modules
has independent interest as
it connects the Emerton-Gee stacks
to objects in geometric representation theory
such as the affine Grassmannians.
\end{abstract}

\maketitle

\tableofcontents

\section{Introduction}

Let $F$ be an arbitrary $p$-adic field,
and let $G$ be an arbitrary reductive group
over $F$ that splits over $L$.
Set $\lsup LG:=\bG\rtimes\Gal(L/F)$.
We establish the following theorem:

\begin{thm}[Corollary \ref{cor:representable},
Theorem \ref{thm:Sha},
Theorem \ref{thm:rep-parabolic}]
\label{thm:main}
(1)
The morphism of Emerton-Gee stacks
$\cX_{F, \lsup LG}\to\cX_{F,\GL_N}$
is relatively representable
by algebraic stacks of finite presentation over $\spf\Z_p$
for any embedding $\lsup LG\to\GL_N$.

(2)
Let $F/E$ be a finite extension.
There is a tautological isomorphism
$\cX_{F, \lsup LG}\cong\cX_{E, \lsup L\Res_{F/E}G}$.

(3)
Let $\lsup LP=\bP\rtimes \Gal(L/F)\subset\lsup LG$
be a parabolic subgroup.
Then
$\cX_{F, \lsup LP}\cong \varinjlim_m
\cX_{F, \lsup LP}^{m}$
where each
$\cX_{F, \lsup LP}^{m}$
is an algebraic stack of finite presentation
over $\spf\Z_p$
with transition maps being nilpotent thickenings.

(4)
For an arbitrary uniformizer of $F$
defining
the moduli stack
of Kummer \'etale $\varphi$-modules
$\clR_{F,\lsup LG}^{\Kummer}$,
the morphism of
$\clR_{F,\lsup LG}^{\Kummer}\to\clR_{F,\GL_N}^{\Kummer}$
is relatively representable
by algebraic stacks of finite presentation over $\spf\Z_p$
for any embedding $\lsup LG\to\GL_N$.
\end{thm}

\begin{remark}
(1)
In the previous version of this manuscript,
we proved the same results 
for tamely ramified groups,
using the same method.
The author later realized that the method
is already strong enough to cover
all reductive groups,
in a conversation with Min.

(2)
An algebraic stack over $\spf \Z_p$
is by definition an algebraic stack
over $\spec \Z/p^a$
for some integer $a>0$.

(3)
In \cite{Min25},
using a completely different strategy,
Min establishes
that
$\cX_{\lsup LG}\to\cX_{\GL_N}$
is relatively representable
by limit-preserving,
locally Noetherian,
formal algebraic stacks over $\spf\Z_p$.
We note that
\cite{Min25} is using the vocabulary of derived
formal algebraic stacks,
where
``locally of finite presentation''
is by definition
``limit-preserving''
as the classic notion of finite presentation
does not make sense in the derived world.

(4)
As is noted by \cite{Min25},
it is hard to obtain (Ind-)representability
for $\clR_{F,\lsup LG}$
using their methods,
as it requires a study of height theory.
This is crucial for constructing the 
potentially semistable substacks of $\cX_{F, \lsup LG}$.

(5)
Our methods are exclusive to reductive groups,
while \cite{Min25} works for much more general groups.

(6)
Our methods and that of \cite{Min25}
are completely orthogonal,
as we don't need intermediate lemmas from each other.
The general principle is that
to get a representability result,
you need some form of finiteness input.
\cite{Min25} uses the $\Gamma$-action to obtain finiteness,
while we use height theory.
At this point, both approaches
have their own serious drawbacks
and complement each other.
\end{remark}

The main difficulty for establishing Theorem \ref{thm:main}
is that we are forced to directly analyze
moduli of $\varphi$-modules
without $\varphi$-stable integral models.
Let $\varphi\in \End(\Z_q/p^a(\!(u)\!))$
such that $\varphi(u)=u^p$ mod $p$,
such that $\varphi(u)\not\in \Z_q/p^a[\![u]\!]$.
The methods of \cite{PR09} and \cite{EG23}
break down
because they require the existence
and Ind-properness
of the Kisin resolution $\cC_{\lsup LG}\to\clR_{\lsup LG}$
of the moduli of
(generalized) \'etale $\varphi$-modules
$\clR_{\lsup LG}$
to establish its representability.
Specializing to either cyclotomic or Kummer \'etale $\varphi$-modules,
when $G$ is a ramified reductive group,
then the na\"ive $\cC_{\lsup LG}$ is always empty:

\begin{center}
\begin{tabular}{lccc}
\toprule
\textbf{Kisin resolution} & \textbf{Unramified} & \textbf{Tamely ramified} & \textbf{Wildly ramified} \\
\midrule
{\textbf{Cyclotomic}}
  & Wach module           & N/A                   & N/A \\
\midrule
{\textbf{Kummer}}
    & Breuil-Kisin module   & parahoric Kisin module & N/A \\
\bottomrule
\end{tabular}
\end{center}
In the tame Kummer setup,
it is possible to fix it by replacing $\cC_{\lsup LG}$
by the moduli of parahoric Kisin modules,
but no such a workaround is known previously
in wildly ramified situations.

We overcome this difficulty
by establishing a framework
of {\it height theory}
(see Definition \ref{def:height})
that allows us to do reductions
in {\it both directions}
--
both
going up and
down {a finite \'etale cover}.
In particular,
we are able to define the closed substacks
of finite height,
and prove

\begin{thm}
[Corollary \ref{thm:RG}]
If $\varphi$
is a generalized Frobenius admitting a height theory,
then
$\calR_{\lsup LG}^{\le h}$ is an algebraic stack
of finite presentation over $\spec \Z/p^a$,
and there exists a $2$-equivalence
\[
\varinjlim_h \clR_{\lsup LG}^{\le h}
\cong \clR_{\lsup LG}.
\]
\end{thm}

Our theorem applies to both cyclotomic $\varphi$-modules
and Kummer $\varphi$-modules
because our framework allows us to {\it go down}
to the case of $F=\Qp$
where $\Z_p[\![u]\!]$ is indeed $\varphi$-stable:

\begin{thm}[Theorem \ref{thm:propagate}]
Suppose $\Lambda$ is a finite local $\Z_p/p^a$-algebra.
Let $\Lambda(\!(u_1)\!)\hookrightarrow \Lambda'(\!(u_2)\!)$
be a finite \'etale extension.
Fix a Frobenius $\varphi_1$ acting on $\Lambda(\!(u_1)\!)$.
Then $\varphi_1$
extends uniquely to $\varphi_2\in \End(\Lambda'(\!(u_2)\!))$,
and $\varphi_1$ admits a height theory
if and only if $\varphi_2$ admits
a height theory.
\label{thm:propagate}
\end{thm}

The representability of $\clR_{\lsup LG}$
is important
because it receives morphisms
from modifications of the affine flag varieties/Grassmannians,
which
$\cX_{\lsup LG}$ does not,
and thus plays a vital role in \cite{LLHLM}
and many follow-up works.

\subsection{Acknowledgement}
We thank Yu Min for very helpful discussions.
We also thank
Matthew Emerton,
Toby Gee, 
Bao Le Hung, 
Stefano Morra, and
 David Savitt
for correspondences and useful discussions
at various stages of this work.

\section{Generalized $\varphi$-modules: examples
and theory}

\begin{exam}[Breuil-Kisin modules]
\label{ex:BK}
Let $F/\Q_p$ be a finite extension.
Let $\varpi$ be a uniformizer of $F$
and let $\varpi^\flat=(\varpi^{1/p^n})_n$
be a compatible choice
of $p$-power roots of $\varpi$.
Set $F_\infty:=F(\varpi^{1/p^\infty})$.
Recall that
$\A_\Inf=W(\cO_{\C^\flat})$
and
$\Z_q[\![ [\varpi^\flat] ]\!]\hookrightarrow
\A_\Inf^{\Gal_{F_\infty}}$
where $q=p^f$ is the size of the residue field of $F$.
Set $v:=[\varpi^\flat]$
and $\cO_\cE:= \Z_q(\!(v)\!)^{\wedge p}\subset \A_\Inf$.

Then $\varphi\in \End(\cO_\cE)$ acts on $\Z_q$ 
as the usual Frobenius,
and sends $v$ to $v^p$.
In particular, $\fS_F:=\Z_q[\![v]\!]$
is $\varphi$-stable.
\end{exam}

\begin{exam}[Cyclotomic $(\varphi, \Gamma)$-modules]
\label{ex:cyc}
Let $F(\zeta_\infty)$
be the union of $p$-power cyclotomic extensions
of $F$
and write $\kappa_{F, \infty}$
for its residue field.

We have the imperfect field of norms
$\E_{F}'\cong \kappa_{F, \infty}(\!(T)\!)$,
which admits a thickening
$\A_F'\cong W(\kappa_{F,\infty})(\!(T)\!)$.

If $F=\Q_p$, then
$\A_{\Qp}'=\Z_p(\!( [\zeta_\infty^\flat]-1 )\!)
\subset W(\C^\flat)$
and
we simply set $T = [\zeta_\infty^\flat]-1$;
as
$\varphi$ acts by $\varphi(T) = (T+1)^p-1$,
$\A_{\Q_p}^{\prime+}:=\Z_p[\![T]\!]$
is $\varphi$-stable.

If $F=\Q_{p^f}[(-p)^{1/e}]$
is tamely ramified over $\Q_p$,
where $e=p^f-1$.
Then, by Krasner's lemma, we 
have 
$\A_{F}'=\Z_{p^f}(\!( ([\zeta_\infty^\flat]-1)^{1/e} )\!)
\subset W(\C^\flat)$;
so, we set
$T = ([\zeta_\infty^\flat]-1)^{1/e}$,
and we have
\begin{align*}
&\varphi(T^e)= (T^e+1)^p-1
\\
\Rightarrow
&\varphi(T) = ((T^e+1)^p-1)^{1/e}
=T^p(1+\frac{p}{T^e}+\dots+\frac{p}{T^{e(p-1)}})^{1/e}
\end{align*}
In particular,
$\A_{F}^{\prime+}:=\Z_{p^f}[\![T]\!]$
is not $\varphi$-stable.

When $F/\Q_p$ is wildly ramified,
we no longer have a simple description
of the $\varphi$-structure.

We will fix a separable closure
$\E^{\sep}/\E_F'$ once for all.
We also fix an unramified extension
$\A^{\sep}/\A_F'$
whose reduction mod $p$
is $\E^\sep$.
Let $\Delta_F$ be the torsion subgroup of
$\Gal(F(\zeta_\infty)/F)$.
We set $\E_F:=(\E_F')^{\Delta_F}$
and $\A_F:=(\A_F')^{\Delta_F}$.
\end{exam}

\begin{exam}[A wildly ramified
extension of Breuil-Kisin $\fS$]
Now we exhibit what happens when $F/\Q_p$
is wildly ramified,
but in the more tractable situation of Breuil-Kisin modules.

In this example, we set $p=2$
and we keep notations as in Example \ref{ex:BK}.

Let $\Omega:= \F_q(\!(v)\!)[X]/(X^p-X-\frac{1}{v})$
be an Artin-Schreier extension.
Set $u:=\frac{1}{X}$
and we have
\begin{align*}
&u^2+v u - v=0.
\\
\Leftrightarrow
&
v=\frac{u^2}{1-u}
=u^2(1+u+u^2+\dots).
\end{align*}
Let $\{1,\gamma\}$ be the Galois group
of $\Omega/\F_q(\!(v)\!)$,
and we have 
\[
\gamma(X)=X+1
\Leftrightarrow
\gamma(u)=\frac{u}{1-u}
\]
and $u \gamma(u)=v$.
Next, we describe an extension
$\wt \Omega/\fS$:
define
\[
\wt\Omega:=\fS[\![u]\!]/(-u^2(1+u+\dots)-v)
\]
where the negative sign is chosen intentionally!
Define
\[
\gamma(u)=\frac{-u}{1-u}
\]
in $\wt \Omega$.
Indeed, the negative sign is crucial
to guarantee $\gamma(\gamma(u))=u$
in $\wt \Omega$.
The reader can check that
$\wt\Omega$
is finite free of rank $2$ over $\fS$.
In particular,
$(\wt \Omega[\frac{1}{v}])^{\wedge p}$
is an unramified
extension of $\cO_\cE$
with residue field $\Omega$
(which is unique up to unique isomorphism
by Hensel's lemma).
By Hensel's lemma,
the $\varphi$-structure on $\Omega$
lifts uniquely to $(\wt \Omega[\frac{1}{v}])^{\wedge p}$.
From $u^2-uv+v=0$,
we have
\[
\varphi(u)=\frac{v^p\pm\sqrt{v^{2p}-4v^p}}{2}
=\frac{u^4}{(1-u)^2}(1\pm \sqrt{1-\frac{4}{u^2}+\frac{8}{u^3}-\frac{4}{u^4}})/2.
\]
in $(\wt \Omega[\frac{1}{v}])^{\wedge p}$
and $\wt \Omega$ is not $\varphi$-stable.
\end{exam}

All examples above play an important role
in the theory of Galois representations.
We list some common features of these examples:

\begin{dfn} [Pivot element]
\label{def:pivot}
Let $\Lambda$ be a $\Z_p$-algebra
with $\varpi\in \Lambda$ such that $p|\varpi^m$
for some $m>0$.
Let $\varphi:\Lambda(\!(u)\!)\to \Lambda(\!(u)\!)$ be a
ring endomorphism
such that $\varphi(u)\in u^{q}+ \varpi\Lambda(\!(u)\!)$
for some $q=p^f$.

A {\bf pivot element} for $\varphi$
is an element $v\in \Lambda(\!(u)\!)$
such that
\begin{itemize}
\item[(H0a)] $v$ is invertible in $\Lambda(\!(u)\!)$;
\item[(H0b)] $\Lambda(\!(u)\!)$ is $v$-adically complete;
\item[(H1)] $\Lambda(\!(u)\!)$ is a finite \'etale
$\Lambda(\!(v)\!)$-module;
\item[(H2)] $\varphi(\Lambda[\![v]\!])\subset \Lambda[\![v]\!]$;
\item[(H3)] There exists an element $u'\in \Lambda(\!(u)\!)$
such that $\Lambda(\!(u')\!)=\Lambda(\!(u)\!)$
and $\Lambda[\![v]\!]\subset \Lambda[\![u']\!]$;
\item[(H4)]
$\Lambda(\!(v)\!)\otimes_{\varphi, \Lambda(\!(v)\!)}
\Lambda(\!(u)\!)
\xrightarrow{x\otimes y\mapsto x\varphi(y)}\Lambda(\!(u)\!)$
is a bijection.
\end{itemize}
Note that (H0a) and (H0b) guarantee that
the inclusion $\Lambda[v]\subset \Lambda(\!(u)\!)$
factors through $\Lambda(\!(v)\!)$.
\end{dfn}

\begin{dfn} [Height theory]
\label{def:height}
We say $\varphi\in\End(\Lambda(\!(u)\!))$
admits a {\bf height theory}
if there exists a finite \'etale extension
$\Lambda(\!(u)\!)\hookrightarrow \wt\Lambda(\!(\wt u)\!)$
such that
$\Lambda(\!(\wt u)\!)$
admits a $\varphi$-structure
with a pivot element $v$.
\end{dfn}

The following theorem
propagates height theories
along finite \'etale extension.

\begin{thm}
Suppose $\Lambda$ is a finite $\Z_p/p^a$-algebra
with residue field $\F$.
Let $\Lambda(\!(u_1)\!)\hookrightarrow \Lambda'(\!(u_2)\!)$
be a finite \'etale extension.
Fix a Frobenius $\varphi_1$ acting on $\Lambda(\!(u_1)\!)$.
Then $\varphi_1$
extends uniquely to $\varphi_2\in \End(\Lambda'(\!(u_2)\!))$,
and $\varphi_1$ admits a height theory
if and only if $\varphi_2$ admits
a height theory.
\label{thm:propagate}
\end{thm}

\begin{proof}
One direction is obvious.
Assume $\varphi_1$ admits a height theory.
Without loss of generality,
we can assume $\Lambda(\!(u_1)\!)$
already admits a pivot element $v$,
and we can assume $v\in \Lambda[\![u_1]\!]$.
It is also harmless to enlarge $\Lambda'(\!(u_2)\!)$
and assume
$\F'(\!(u_2)\!)/\F(\!(u_1)\!)$
is a finite Galois extension,
and by induction on its extension degree
(the absolute Galois group
of a local field in char $p$ is solvable),
we can assume it is a cyclic Galois extension
of prime degree $r$.
The only hard case is when $r=p$
and $\Lambda'=\Lambda$,
and we assume both.
An irreducible Artin-Schreier equation
over $\F(\!(u_1)\!)$
is of the form
\[
X^p-X=\frac{a(u_1)}{u_1^n}
\]
where the constant term of $a(u_1)\in
\F[\![u_1]\!]$ is non-zero.
We lift $a(u_1)$
to an invertible element
$\wt a(u_1)\in \Lambda[\![u_1]\!]$,
and consider the equation
\[
g(U):=U^p + u_1^n\wt a(u_1)^{-1} U^{p-1}- u_1^n\wt a(u_1)^{-1}
\in \Lambda[\![u_1]\!][U].
\]
It is clear that
\[
\Lambda[\![U]\!]:=\Lambda[\![u_1]\!][U]/(g(U))
\]
is finite free of rank $p$ over 
$\Lambda[\![u_1]\!]$.
Note that
\[
\Lambda[\![U]\!][1/U] \equiv \F(\!(u_2)\!)
\mod m_\Lambda,
\]
and thus
$\Lambda[\![U]\!][1/U]
\cong \Lambda(\!(u_2)\!)$
as $\Lambda(\!(u_1)\!)$-algebra
by Hensel's lemma
(or, the topological invariance of the finite \'etale site).
Since $\varphi(g(U))$ mod $m_\Lambda$
(as a polynomial in variable $\varphi(U)$)
has no repeated
roots, and thus by Hensel's lemma
$u_2^q$ lifts uniquely
to a solution of $\varphi(g(U))=0$.
Finally,
$v\in \Lambda[\![u_1]\!]\subset \Lambda[\![U]\!]$
is also a pivot element of for $\varphi_2$.
\end{proof}

\subsection{Moduli of \'etale $\varphi$-modules}
Let $H$ be a smooth affine group scheme over $\Z_p$.
An \'etale $\varphi$-module with $H$-structure
with $R$-coefficients
is an $H$-torsor $T$ over $\Lambda\otimes_{\Z_p}R(\!(u)\!)$
together with an isomorphism
$\varphi_T:\varphi^*T\xrightarrow{\cong} T$.

\begin{dfn}
Denote by $\clR_H^\varphi$
for the moduli stack of \'etale $\varphi$-modules
with $H$-structure,
defined over $\spf \Z_p$.
\end{dfn}

The following theorem shows
that to prove representability theorems,
it is harmless
to replace $\Lambda(\!(u)\!)$
by a finite unramified extension.

\begin{thm}
Keep notations as in Theorem \ref{thm:propagate}.
If $\cZ\subset \clR_{H}^{\varphi_2}$
is an algebraic closed substack
of finite presentation over $\spec \Z/p^a$,
then
$\clR_H^{\varphi_1}
\underset{\clR_H^{\varphi_2}}{\times}
\cZ$
is an algebraic stack of finite presentation over
$\spec \Z/p^a$.
\label{thm:restrict}
\end{thm}

\begin{proof}
The proof is the same as \cite[Lemma 3.7.5]{EG23}.
The idea is that $\clR_H^{\varphi_1}
\underset{\clR_H^{\varphi_2}}{\times}\cZ$
can be presented as
a fiber product involving
only $\cZ$
and morphisms between (fiber products)
of $\cZ$.
\end{proof}

\subsection{Langlands parameters
and non-split groups}

Let $G$ be a reductive group over $F$
that splits over $L$.
We have $\lsup LG=\bG\rtimes\Gal(L/F)$
where $\bG$ is a split connected reductive group
whose root system is dual to that of $G$.
Recall that a Langlands parameter for $G$
is a group homomorphism
$\Gal_F\to \lsup LG(\bC)$
such that the composite
$\Gal_F\to \lsup LG(\bC)\to \Gal(L/F)$
is the quotient map.

Now we go back to consider \'etale $\varphi$-modules
with $H$-structure
where $H=H^\circ\rtimes \Delta$
is a semidirect product of a connected reductive group
with a finite group.
It is natural to have the following definition:

\begin{dfn}
Fix one distinguished element
$\std: \spf \Z_p\to \clR_\Delta^\varphi$
once for all.
Define
\[
\clR_H=\clR_{H^\circ\rtimes \Delta}:=
\clR^\varphi_H \times_{\clR^\varphi_\Delta, \std}\spf \Z_p.
\]
In particular,
if $H$ is already connected,
then $\clR_H=\clR_H^\varphi$.
We write $\clR_H={\clR_H}_{/\Lambda(\!(u)\!)}$
to emphasize the dependence on 
$\Lambda(\!(u)\!)$.

The distinguished element
$\std:\spf \Z_p\to \clR_\Delta^\varphi$
defines
a $\Delta$-torsor
$T_{\std}$
over $\Lambda(\!(u)\!)$.
Indeed,
$T_{\std}=\spec \Lambda'(\!(u')\!)$
is a finite \'etale cover of $\spec \Lambda(\!(u)\!)$
and
$\varphi_{T_{\std}}:
\Lambda'(\!(u')\!)\to \Lambda'(\!(u')\!)$
is an endomorphism extending
$\varphi$.
\end{dfn}

The following proposition
allows us to reduce representability questions
to the case of split groups.

\begin{prop}
If $\cZ\subset {\clR_{H^\circ}}_{/\Lambda'(\!(u')\!)}$
is an algebraic closed substack
of finite presentation over $\spec \Z/p^a$,
then
${\clR_H}_{/\Lambda(\!(u)\!)}
\underset{{\clR_{H^\circ}}_{/\Lambda'(\!(u')\!)}}{\times}
\cZ$
is an algebraic stack of finite presentation over
$\spec \Z/p^a$.
\label{prop:restrict}
\end{prop}

\begin{proof}
So,
$T_{\std}\otimes_\Lambda \Lambda'
=\spec \Lambda'\otimes_\Lambda\Lambda'(\!(u')\!)$
is a trivial $\Delta$-torsor
over $\spec \Lambda'(\!(u')\!)$.

We fix a marked point
$\Pt:\spec \Lambda'(\!(u')\!)\to T_{\std}\otimes_\Lambda \Lambda
=\prod \spec \Lambda'(\!(u')\!)$
once for all.
For each object $(T, \varphi_T)$
of ${\clR_{H}}_{/\Lambda(\!(u)\!)}(A)$
where $A$ is a $\Lambda'$-algebra,
there is a tautological morphism
\[
T\to 
T\times^{\lsup LH}\Delta=T_{\std}\otimes_{\Z_p(\!(u)\!)}A(\!(u)\!)
=\spec \Lambda'\otimes A(\!(u')\!)
=\spec \Lambda'\otimes_\Lambda(A\otimes_{\Z_p}\Lambda)(\!(u)\!).
\]
Pulling back of the marked point 
of $T\times^{\lsup LH}\Delta$
defines
a $H^\circ$-torsor
$T^\circ\subset T$
over $(A\otimes_{\Z_p}\Lambda)(\!(u')\!)$.
Moreover, there is an identification
\[
T^\circ \otimes_\Lambda\Lambda'
=T.
\]
As a consequence,
$T$ can be regarded
as a $H^\circ$-torsor
over $(A\otimes_{\Z_p}\Lambda')(\!(u')\!)$,
and is tautologically an object
of 
$\clR^{\varphi_{u'}}_{H^\circ}(A)$.

Therefore we have explicitly defined an isomorphism
\[
{\clR_H}_{/\Lambda'(\!(u')\!)}
\xrightarrow{\cong}
\clR^{\varphi_{u'}}_{H^\circ}
=
{\clR_{H^\circ}}_{/\Lambda'(\!(u')\!)}.
\]
The rest of the proof is the same as that for
Theorem \ref{thm:restrict}.
\end{proof}

\section{Frobenius stabilization of loop groups}

In this section,
fix a Frobenius endomorphism
$\varphi\in \End(\Z_q/p^a(\!(u)\!))$
such that $\varphi(u)=u^p+O(p)$.
Let $H$ be a connected reductive group scheme over $\spec \Z_p$
and set $\clH:=\Res_{\Z_q/\Z_p}H$.
Fix an embedding $i:\clH \hookrightarrow\GL_N$.
We assume $\varphi$ admits a pivot element $v
\in \Z_q/p^a[\![u]\!]$.

\begin{dfn}[Principle congruence subgroups
and bounded height subfunctors]
The positive loop group 
\[
L^+\clH: A\mapsto\clH(A[\![u]\!])
\] has
congruence subgroups $U_n$,
defined by
$$
U_n: A \mapsto \{X\in L^+\clH(A)| i(X) \in 1+ v^{n}\operatorname{Mat}_{N\times N}(A[\![u]\!]) \}.
$$
Also define
subfunctors $L^{\le h}\clH$ of the loop group
$L\clH:A\mapsto\clH(A(\!(u)\!))$
by
\[
L^{\le h}\clH: 
A\mapsto\{X\in L\clH(A)| i(X),i(X)^{-1}\in v^{-h}\operatorname{Mat}_{N\times N}(A[\![u]\!]) \}
\]
\end{dfn}

We warn the reader that $\varphi(L^+\clH)\not\subset L^+\clH$
in general.
We define $\varphi$-twisted conjugation of
$L\clH$ on itself
by $(x,y)\mapsto \varphi(y) x y^{-1}$.

\begin{lem}
(1)
If $n\in \Z$ is sufficiently large, then
$\varphi(U_n)\subset L^+\clH$.
In particular,
the quotient stack
$[L^{\le h}\clH/_\varphi U_n]$
is well-defined.

(2)
If $n\in \Z$ is sufficiently large, then
$
[L^{\le h}\clH/U_n] \to [L^{\le h}\clH/_\varphi U_n]
$
sending $X\mapsto X$
is an equivalence of categories.
($U_n$ acts by left-translation if notation is not decorated.)
\end{lem}

\begin{proof}
(1) is obvious.
(2) follows from a standard argument by Pappas-Rapoport.
See the proof of \cite[Prop. 2.2]{PR09} or \cite[Lemma 5.2.9]{EG23}.
\end{proof}

The factor group $\clH_n:= L^+\clH/U_n \cong
\Res_{(\Z_p[\![u]\!]/u^{n'})/\Z_p}(\clH\otimes 
\Z_p[\![u]\!]/u^{n'})$
is smooth since it is the Weil restriction of a smooth group;
we note that Weil restrictions
preserve smoothness, but do not preserve properness
in general
(unless along a finite \'etale cover).

\begin{cor}
\label{cor:PR0}
If $n$ is sufficiently large, the quotient stack
$$
[L\clH/_\varphi U_n] = \dirlim{h}~[L^{\le h}\clH/_\varphi U_n]
$$
is an Ind-algebraic stack with closed
algebraic substacks $[L^{\le h}\clH/_\varphi U_n]$
that are finitely presented over $\Z/p^a$.
\end{cor}

\begin{proof}
$[L^{\le h}\clH/ U_n]$
is a $\clH_n$-torsor over
$[L^{\le h}\clH/ L^+\clH]$,
which is a closed subscheme of the affine Grassmannian.
\end{proof}

Next we consider the non-standard filtration on $[L\clH/_\varphi U_n]$.

\begin{dfn}
\label{def:frob-stab-LmG}
There exists an integer $c_\varphi > 0$
such that
$\varphi(\Z_q/p^a[\![u]\!])\subset
\frac{1}{v^{c_\varphi}}\Z_q/p^a[\![u]\!]$.
Define the stabilized filtration as follows:
\begin{align*}
L^{\le h}_\varphi\clH: A\mapsto &\{X\in L\clH(A)|p^{a-i} i(X), p^{a-i} i(X)^{-1}\in v^{-h - c_\varphi i}\operatorname{Mat}_{N\times N}(\Lambda[\![u]\!]),~ 0\le i\le a\}.
\end{align*}
\end{dfn}

\begin{lem}
\label{lem:pl-1}
(1)
$L^{\le h}_\varphi\clH$ is a closed subscheme of $L\clH$.

(2)
$L^{\le h}\clH\subset L^{\le h}_\varphi\clH \subset L^{\le h+ a c_\varphi}\clH$.

(3)
$L^{\le h}_\varphi\clH$ is closed under the $\varphi$-twisted
conjugation action of $L^+\clH$.
\end{lem}

\begin{proof}
(1):
Let $i_s:L\clH \to L \operatorname{Mat}_{n\times n}
\xrightarrow{x \mapsto p^s x} L \operatorname{Mat}_{n\times n}$
be the multiplication by $p^s$ morphism.
Note that
$L^{\le h}_\varphi\clH$ is the intersection
of the pullbacks of
closed subschemes
along $i_s$.

(2):
Immediate from the definition.

(3): Note that $\varphi(L^+\clH(A))\subset L^+\clH(A) 
+ \frac{p}{v^{c_\varphi}}\operatorname{Mat}_{N\times N}(A[\![u]\!])$.
Let $g\in L^+\clH(A)$
and $x\in L^{\le h}_\varphi\clH(A)$.
Write $\varphi(g)=g_1 + \frac{p}{v^{c_\varphi}}g_2$
where $g_1, g_2\in \operatorname{Mat}_{N\times N}(A[\![u]\!])$.
Note that
$$
g^{-1} p^{a-i} x \varphi(g)
=g^{-1} p^{a-i} x g_1 + g^{-1} \frac{p^{a-i+1} x}{v^{c_\varphi}} g_2
\in v^{-h - c_\varphi i}\operatorname{Mat}_{N\times N}(A[\![u]\!]),
\text{~and}
$$
$$
g^{-1} p^{a-i} x^{-1} \varphi(g)
=g^{-1} p^{a-i} x^{-1} g_1 + g^{-1} \frac{p^{a-i+1} x^{-1}}{v^{c_\varphi}} g_2
\in v^{-h - c_\varphi i}\operatorname{Mat}_{N\times N}(A[\![u]\!]).
$$
So we are done.
\end{proof}

\begin{thm}
\label{thm:PR-0}
The stack $[L\clH/_\varphi L^+\clH]$
is representable by an Ind-algebraic stack.
We have an Ind-presentation
$$
[L\clH/_\varphi L^+\clH]= \dirlim{h}~[L^{\le h}_\varphi\clH/_\varphi L^+\clH]
$$
where each $[L^{\le h}_\varphi\clH/_\varphi L^+\clH]$
is an algebraic stack of finite type over $\Z/p^a$.
\end{thm}

\begin{proof}
It follows from the diagram
$$
\xymatrix{
    [L^{\le h}\clH/U_{n}] \ar[r]^{\cong}
    \ar[d]^{\text{$\clH_{n}$-torsor}} &
    [L^{\le h}\clH/_\varphi U_{n}] \ar@{^{(}->}[d] \\
    [L^{\le h}\clH/L^+\clH] 
    & \dirlim{h}[L^{\le h}\clH/_\varphi U_{n}]
    \ar@{=}[r]& \dirlim{h}[L^{\le h}_\varphi\clH/_\varphi U_{n}]
    &
    [L^{\le h}_\varphi\clH/_\varphi U_{n}]
     \ar@{^{(}->}[l]\ar[d]^{\text{$\clH_n$-torsor}}\\
 &       && [L^{\le h}_\varphi \clH/_\varphi L^+\clH]
}
$$
that each $[L^{\le h}_\varphi \clH/_\varphi L^+\clH]$
is an algebraic stack of finite type over $\Z/p^a$.
\end{proof}

Both $L^{\le h}\clH$
and $L^{\le h}_\varphi\clH$
are purely technical devices
to get the Ind-algebraicity of
$[L\clH/_\varphi L^+\clH]$.

\begin{dfn}
Define the Kisin stack $\cC_{\clH}^{\le h}$
to be the scheme-theoretic image of
$L^{\le h}\clH$
in 
$[L\clH/_\varphi L^+\clH]$,
which is an algebraic stack of finite type over $\Z/p^a$.
We also define
$\L^{\le h}\clH$
as the fiber product
\[
\L^{\le h}\clH
:=
\cC_{\clH}^{\le h}
\underset{[L\clH/_\varphi L^+\clH]}{\times}
L\clH,
\]
which is supposed to be a more correct
filtration than both
$L^{\le h}$
and $L^{\le h}_\varphi$.

There is a natural forgetful morphism
$\pi_h:\cC_{\clH}^{\le h}\to \clR_H$.
Write $\clR_{H}^{\le h}$
for the scheme-theoretic image of $\pi_h$.
\end{dfn}

\begin{thm}
\label{thm:pappas-rapoport}
(1) 
The morphism $\pi_h$
is representable by algebraic spaces, proper, and of finite presentation
over $\Z/p^a$,
and
$\cC^{\le h}_{\clH}$
is an algebraic stack of finite presentation
over $\Z/p^a$ with affine diagonal.

(2)
The diagonal morphism
$\clR_{H} \to \clR_{H}\times\clR_{H}$
is representable by algebraic spaces,
affine, and of finite presentation.
Moreover,
for each morphism $\spec A\to \clR_{H}\times\clR_{H}$,
the change of group morphism
$$
\spec A\underset{\clR_{H}\times\clR_{H}}{\times}\clR_{H}
\to
\spec A\underset{\clR_{\GL_N}\times\clR_{\GL_N}}{\times}\clR_{\GL_N}
$$
is a closed immersion.
\end{thm}

\begin{proof}
(1)
See \cite[Corollary 2.6]{PR09},
and the key is that $\pi_h$
is relatively representable by a closed subscheme
of the $\cF$-twisted affine Grassmannian
(we can't find a reference
with sufficient level of generality,
so we included it in the appendix: Lemma
\ref{gr-lem3}).

(2)
The same proof as \cite[Theorem 5.4.9(3)]{EG21}
works.
\end{proof}

\begin{prop}
(1)
$\clR_{H}$
admits a versal ring at all finite type points.

(2)
$\clR_{H}^{\le h}$
admits a Noetherian versal ring at all finite type points.
\label{prop:versal}
\end{prop}

\begin{proof}
(1)
By Steinberg's theorem,
all $\clH$-torsors over $W(\bar\F_p)/p^a(\!(u)\!)$
are trivial $\clH$-torsors.
So, just use the completion
of $L\clH$
at a finite type point.

(2)
The $n$-th principal congruence subgroup $U_n$
is smooth because it is obtained
as the dilation of
smooth groups.
Note that
for $n\gg h$,
$[\L^{\le h}\clH/_\varphi U_n]
=\L^{\le h}\clH/U_n$.

Let $x\in\clR^{\le h}_{H}(\bar\F_p)$
be a finite type point,
and let $R_x$ be a (non-Noetherian)
versal ring of $\clR_{H}$
constructed in part (1).
Write $R_x^\cC$
for the scheme-theoretic image
of $\spf R_x\times_{\clR_{H}}\cC^{\le h}_{\clH}$
in $\clR_{H}$.
By the universal property
of scheme-theoretic image,
$U_n$ acts on $\spf R_x^{\cC}$,
and $\spf R_x^\cC/U_n\to \clR_{H}$
factors through
a Noetherian versal ring of $\clR^{\le h}_{H}$
at $x$
by \cite[Lemma 3.2.16]{EG21}.
Our proof is identical to that in \cite[Theorem 5.4.19]{EG21}.
Also see the previous version of this manuscript
for an elaborated proof
(which is now removed due to readers complaining
it being tedius).
\end{proof}

\subsection{The universal family of Kisin lattices
over $\cC^{\le m}_{\clH}$}

From now on,
we base change all stacks from $\Z/p^a$
to $W(\bar\F_p)/p^a$
to simplify exposition.
We also omit subscripts when it is clear from the context.
In particular,
$\clH=\Res_{\Z_q/\Z_p}H=H^{\times[\F_q:\F_p]}=H \times \dots\times H$
is just a product of $H$
after base change.

We describe the universal family in concrete terms.
Set
\[
\cC^{\le m}[\![ u]\!]
:= 
\cC^{\le m}\times \spec \Z_p[\![u]\!],
\]
which is an algebraic stack
(of infinite type, not seen as formal algebrai stack
of topological finite type!).
The universal family
is a pair $(T^\univ, \varphi_T^\univ)$
where $T^\univ$
is a $\clH=H \times \dots\times H$-torsor
over $\cC^{\le m}[\![ u]\!]$,
and
$\varphi_T^\univ:
\varphi^*T^\univ[\frac{1}{u}] \xrightarrow{\cong}T^\univ
[\frac{1}{u}]$
is an isomorphism,
but not (!) an $\clH$-torsor morphism.
To elaborate,
there is an outer automorphism
$\sigma: \clH \to \clH$
that permutes coordinates
$(g_1, \dots, g_f)\mapsto (g_2, \dots, g_f, g_1)$,
and we
write $\sigma: T^\univ\xrightarrow{\cong} T^\univ$
for the morphism
$T^{\univ}\to T^{\univ}\times^{\clH, \sigma}\clH
=T^\univ$.
The composite $\sigma \varphi_T^\univ$
is a $\clH$-torsor isomorphism.

This perspective is a little bit
inconvenient;
alternatively,
equip
$\varphi_\sigma^*T^\univ$($:=
\varphi^*T^\univ$
as an algebraic stack)
with the shifted $\clH$-torsor structure.
Then
$\varphi_T^\univ:
\varphi^*_\sigma T^\univ[\frac{1}{u}] \xrightarrow{\cong}T^\univ
[\frac{1}{u}]$
is a genuine $\clH$-torsor isomorphism.

We remark that, a priori,
the universal family only
exists
as a formal algebraic stack;
and it is the Grothendieck existence theorem
and the Grothendieck algebraization theorem
that allow us to get the algebraic
universal family as above.

\begin{prop}
Let $\spf R\to \clR_{H}^{\le m}$
be a Noetherian versal ring.

The morphism
\[
\widehat{\pi}:\spf R\underset{\clR_{H}}{\times}
\cC^{\le m}
\to \spf R
\]
admits an algebraization
which we denote by
\[
\pi:\cC_R
\to \spec R.
\]
There exists a pair
$(T_R^\cC, \varphi_{T_R^\cC})$
where $T_R^\cC$
is a $\clH$-torsor
over the {\bf scheme}
$\cC_R[\![u]\!]
= \cC_R\times \spec \Z_p[\![u]\!]$
and that
$\varphi_{T_R^\cC}:\varphi^*_\sigma T_R^\cC[1/u]\to T_R^\cC[1/u]$
is a $\clH$-torsor isomorphism.

Define
\[
\fC_R:=\varinjlim_n \cC_{R/m_R^n}(\!(u)\!)
\overset{\wh \j}{\hookrightarrow} \cC_R
,
\]
for the $m_R$-adic completion,
which is a formal scheme
with a projection
$\wh \pi:\fC_R\to\spf R$.

Denote by
$\fT_R^\cC$
the $m_R$-adic completion
of $T_R^\cC$
that lives over $\fC_R$.

Denote by $(\fT_R, \varphi_{\fT_R})$
the \'etale $\varphi$-module
over $\spf R$
that corresponds to the
given morphism $\spf R\to \clR_{H}^{\le m}$.
Then there is an isomorphism
\[
\fT_R^\cC \cong \wh\pi^* \fT_R
\]
respecting $\varphi$-structures.
\label{prop:CR}
\end{prop}

\begin{proof}
We interpret $\clH$-torsors
as tensor functors
using Tannakian formalism,
and note that
the Grothendieck existence theorem
(c.f. \cite[Tag 0CYW]{Stacks})
is an exact functor
(because it is an equivalence of abelian categories)
that preserves tensor products
and duality
for finite projective modules.
The rest is tautological.
\end{proof}

Here, an \'etale $\varphi$-module over $\spf R$
is interpreted
as a compatible sequence
of \'etale $\varphi$-modules
over $\spec R/m_R^n$
(which is standard).

\begin{thm}
Keep notations as in Proposition \ref{prop:CR}.
There exists
an \'etale $\varphi$-module
$(T_R, \varphi_{T_R})$
over $\spec R$
with $\clH$-structure
whose completion is
$(\fT_R, \varphi_{\fT_R})$.

In particular,
all finite type points
of $\clR^{\le m}_{H}$
admit an effective Noetherian versal ring.
\label{thm:effective}
\end{thm}

\begin{proof}
We interpret
$\clH$-torsors as tensor functors.
So, for any algebraic representation
$\clH\to \GL(V)$,
$T_R(V)$
denotes the coherent sheaf
$T_R\times^{\clH}V$.
Write $\wh{R(\!(u)\!)}$
for the $m_R$-adic completion of $R(\!(u)\!)$.
Be careful with the difference between
$\spf \wh{R(\!(u)\!)}$
and $\spec \wh{R(\!(u)\!)}$.

We first forget the $\varphi$-structure,
See Definition \ref{def:fun-ext-vect}
for the definition of 
$\Vect_{\spec \wh {R(\!(u)\!)}, \spec R[\![u]\!], \cC_{R}[\![u]\!]}$.
Define a lax monoidal functor 
\[
(T_R^\cC, \fT_R): \frep_{\clH}\to \Vect_{\spec \wh {R(\!(u)\!)}, \spec R[\![u]\!], \cC_{R}[\![u]\!]}
\]
by sending $V \mapsto (T_R^\cC(V), \fT_R(V))$.
By Proposition \ref{prop:decomp},
there exists a lax monoidal functor
\[
\xi: \Vect_{\spec \wh {R(\!(u)\!)}, \spec R[\![u]\!], 
\cC_{R}[\![u]\!]}\to \Coh_{\spec R(\!(u)\!)}.
\]
Define $T_R:=\xi\circ (T_R^\cC, \fT_R)$
as the composite.

For each algebraic representation $V$, $T_R(V)$ is equipped with an \'etale $\varphi$-structure, as follows.
We first assume 
$\varphi_{T_R^\cC(V)}: \varphi_\sigma^*T_R^\cC(V)[\frac{1}{u}]\to T_R^\cC(V)[\frac{1}{u}]$
is effective in the sense that
$\varphi_{T_R^\cC(V)}(T_R^\cC(V))\subset
T_R^\cC(V)$.
By the proof of Proposition \ref{prop:decomp}
(any unfamiliar notation is defined there),
\[
T_R(V)=j^*\invlim{n} (j, \pi)_*i_n^*(\fT_R(V), T_R^\cC(V))
=j^*\invlim{n}\kappa_n.
\]
Since $T_R^\cC$ admits an effective $\varphi$-structure,
so does $(j, \pi)_*i_n^*(\fT_R(V), T_R^\cC(V))$ for each $n$.
By passing to the the inverse limit and then the localization,
$T_R(V)$ admits a $\varphi$-structure.
In general,
when $T_R^\cC(V)$
is not effective,
we can twist it by a sufficiently high $N$-th Tate twist
to get
effective $T_R^\cC(V)(N)$,
define the $\varphi$-structure,
and transport the $\varphi$-structure
back to $T_R^\cC(V)$
by tensoring with the $(-N)$-th Tate twist.
Here, a Tate twist
can be any finite free $\Lambda(\!(v)\!)$-module
$T(1)$
of rank $1$,
equipped with
$\varphi_{T(1)}:\varphi^*\Lambda(\!(v)\!)\xrightarrow{\cong} \Lambda(\!(v)\!)$
such that $\varphi_{T(1)}(v)\subset v\Lambda[\![v]\!]$.

We still assume the effectivity of $\varphi_{T_R^\cC}$
without loss of generality.
Since the inverse limit functor is left-exact,
$\varphi^* T_R(V)\to T_R(V)$ is injective.
So, $\varphi^* T_R(V)\cong
\varphi(\varphi^* T_R(V))\subset
T_R(V)$ is isomorphic to its image;
here, we omit the subscript from $\varphi$
as there is no ambiguity.
Next, we show its surjectivity
after inverting $u$.
Choose an integer $h>0$ such that
\[
v^hT_R^\cC(V) \subset \varphi(\varphi^*T_R^\cC(V)),
\]
which exists by the following sublemma:

\begin{sublem}
If $\varphi(R[\![v]\!])\subset R[\![v]\!]$
and $\cK$ is a coherent sheaf over
over $X[\![u]\!]$
equipped
with an \'etale $\varphi$-structure
$\varphi:\cK[1/u] \xrightarrow{\cong} \cK[1/u]$
where $X$ is a quasi-compact $R$-scheme,
then there exists an integer $h>0$
such that $v^h\cK \subset \varphi(\varphi^*\cK)$,
where $\varphi(\varphi^*\cK)$ is the injective image of
$R[\![v]\!]\otimes_{\varphi,R[\![v]\!]}\cK
\hookrightarrow
R[\![v]\!]\otimes_{\varphi,R[\![v]\!]}\cK[1/v]
\cong \varphi^*\cK[1/u] \xrightarrow{\cong} \cK[1/u]$.
\end{sublem}

\begin{proof}
Since $X$ is quasi-compact, it suffices to prove the
affine case $X=\spec R$.
Let $\{x_1,\dots, x_n\}$ be a set of generators of $\cK$.
For each $x_i$, there exists an element $v^h$ such that
$v^h\varphi(x_i)\in \cK$, so we are done.
\end{proof}

Since push-forward is left-exact, we have
\[
v^h \pi_*T^\cC_R(V) \subset \pi_*\varphi(\varphi^*T^\cC_R(V))
\cong \pi_*\varphi^*T^\cC_R(V).
\]
By the flat base change theorem and \cite[Lemma 5.2.5]{EG21},
$\pi_*\varphi^*T^\cC_R(V)=\varphi^*\pi_*T^\cC_R(V)$,
and so
\[
v^h \pi_*T^\cC_R(V) \subset \varphi(\varphi^*\pi_*T^\cC_R(V)).
\]
As a consequence $v^h \kappa_n\subset \varphi( \varphi^*\kappa_n)$
for all $n$.
Since the inverse limit is left exact,
\[
v^h\invlim{n} \kappa_n\subset \invlim{n} 
\varphi(\varphi^*\kappa_n);\]
note that $\varphi$ is injective (a bijection
onto its image), and tautologically commutes
with inverse limits:
\[
v^h\invlim{n} \kappa_n\subset \varphi(\invlim{n} 
\varphi^*\kappa_n);\]
inverting $v$, we get
$T_R(V)[\frac{1}{v}]\subset 
\varphi(\varphi^*T_R(V))[\frac{1}{v}]$,
that is, $\varphi_{T_R(V)}:
\varphi^*T_R(V)[1/u]\to T_R(V)[1/u]
$ is an isomorphism.

So far we've shown for each $V$,
$T_R(V)$ is equipped with an \'etale $\varphi$-module
structure,
which is clearly functorial in $V$.
So the functor $T_R$ is upgraded into
a lax monoidal functor from
$\frep_{\clH}$
to the category of
(not necessarily projective) \'etale $\varphi$-modules
over $R(\!(u)\!)$.
Apply the restriction of scalar functor
$\Res_{R(\!(u)\!)}^{R(\!(v)\!)}$
to $T_R$
to get a functor $\Res_{R(\!(u)\!)}^{R(\!(v)\!)}(T_R)$ from
$\frep_{\clH}$
to the category of
\'etale $\varphi$-modules
over $R(\!(v)\!)$.
Note that the restriction of scalar functor
is exact and commutes with completion, but is not lax monoidal;
in particular, the $m_R$-adic completion
of $\Res_{R(\!(u)\!)}^{R(\!(v)\!)}(T_R)$
is an exact functor.
By Proposition \ref{prop:decomp},
the $I$-adic completion of $\Res_{R(\!(u)\!)}^{R(\!(v)\!)}(T_R)$
is $\Res_{\wh {R(\!(u)\!)}}^{\wh {R(\!(v)\!)}}\fT_R$.
So, by Lemma \ref{lem-eff3} (the special case $0\to M=M\to 0\to 0$),
each $\Res_{R(\!(u)\!)}^{R(\!(v)\!)}T_R(V)$ is a finite projective module over $R(\!(v)\!)$.
So, $T_R(V)$ is a finite projective module over $R(\!(u)\!)$,
by noting that if $A\to B$ is finite \'etale, then
being $A$-projective is equivalent to being $B$-projective
as $B$ is a direct summand of $B\otimes_AB$.
By Lemma \ref{lem-eff4}, $T_R$
is a faithful, exact, strict monoidal functor,
and we are done.
\end{proof}

\begin{cor}
\label{thm:RG}
$\calR_{H}^{\le h}$ is an algebraic stack
of finite presentation over $\spec \Z/p^a$,
and there exists a $2$-equivalence
\[
\varinjlim_h \clR_H^{\le h}
\cong \clR_H.
\]
\end{cor}

\begin{proof}
It follows from \cite[Theorem 1.1.1]{EG21}
(whose assumptions are checked in
Theorem \ref{thm:pappas-rapoport},
Proposition \ref{prop:versal}
Theorem \ref{thm:effective}),
except that
\cite[Theorem 1.1.1]{EG21}
only implies $\mathcal{R}_{H}^{\le h}$
is an algebraic stack locally of
finite presentation.
It remains to show
$\mathcal{R}_{H}^{\le h}$
is quasi-compact and quasi-separated.
By the definition of scheme-theoretic image,
$\cC_{\clH}^{\le h}\to \clR_{H}^{\le h}$
is scheme-theoretically dominant.
Since $\cC_{\clH}^{\le h}$ is of finite presentation,
we can choose a finitely presented smooth cover
$\spec S\to \cC_{\clH}^{\le h}$.
Note that the composition $\spec S\to \clR_{H}^{\le h}$
is scheme-theoretically dominant.
So, $\clR_H^{\le h}$ is quasi-compact.
By Theorem \ref{thm:pappas-rapoport}, $\mathcal{R}_{H}^{\le h}$
has affine diagonal and thus is
quasi-separated, as required.
\end{proof}

\begin{cor}
\label{cor:cartesian-RG}
The following diagram is commutative
$$
\xymatrix{
    \cC_{\clH}^{\le h}\ar[d]\ar[r] & \cC_{\GL_N}^{\le h}\ar[d]\\
    \calR_{H}^{\le h}\ar[d]\ar[r] & \calR_{\GL_N}^{\le h}\ar[d]\\
    \calR_{H}\ar[r] & \calR_{\GL_N}
}
$$
where the bottom square is Cartesian,
and 
$$
\cC_{\clH}^{\le h} \to \cC_{\GL_N}^{\le h}\underset{\clR_{\GL_N}^{\le h}}{\times}\clR_{H}^{\le h}
$$
is a closed immersion.
\end{cor}

\begin{proof}
The morphism $\cC_{\clH}^{\le h}\to\calR_{\GL_N}$
factors through $\calR_{H}^{\le h}\times_{\calR_{\GL_N}}\calR_{\GL_N}^{\le h}$,
which is a closed substack of $\calR_{H}$.
By \cite[Lemma 3.2.31]{EG21},
$\calR_{H}^{\le h}$ is a closed substack of
$\calR_{H}^{\le h}\times_{\calR_{\GL_N}}\calR_{\GL_N}^{\le h}$,
which forces
$\calR_{H}^{\le h} = \calR_{H}^{\le h}\times_{\calR_{\GL_N}}\calR_{\GL_N}^{\le h}$.
To show
$\cC_{\clH}^{\le h} \to \cC_{\GL_N}^{\le h}\underset{\clR_{\GL_N}}{\times}\clR_{H}^{\le h}$
is a closed immersion,
we argue by descent and show that
for a chosen smooth cover $\spec B\to \clR_{H}^{\le h}$,
the base-changed morphism
is a closed immersion;
since both schemes above are closed subschemes of the twisted affine Grassmannian
(for $\clH$ and $\GL_N$, respectively) by the proof of
part (1) of Theorem \ref{thm:pappas-rapoport},
the claim follows from Lemma \ref{gr-lem2}.
%
\end{proof}

\begin{cor}
The morphism
$\clR_H\to \clR_{\GL_N}$
is relatively representable by
algebraic stacks of finite presentation
over $\spec \Z/p^a$.
\label{cor:finite}
\end{cor}

\begin{proof}
Since $\clR_{\GL_N}$
is limit preserving,
it suffices to check on finitely presented
test objects
$\spec A\to \clR_{\GL_N}$,
which necessarily factors through some
$\clR_{\GL_N}^{\le h}$.
The corollary now follows from Corollary \ref{cor:cartesian-RG}.
\end{proof}

\section{Representability
for $\clR_{\lsup LG}$
and $\cX_{\lsup LG}$}

\begin{thm}
Let $\varphi\in \Z_q/p^a(\!(u)\!)$
be a Frobenius endomorphism with a height theory
(c.f. Definition \ref{def:height}),
and let $G$ be an arbitrary reductive group
over an arbitrary $p$-adic field $F$.
Then $\clR_{\lsup LG}\to\clR_{\GL_d}$
is relatively representable
by algebraic stacks of finite presentation over $\spec \Z/p^a$
for any embedding $\lsup LG\to\GL_d$.
\end{thm}

\begin{proof}
Combine Theorem 
\ref{thm:propagate},
Theorem \ref{thm:restrict},
and Proposition \ref{prop:restrict}
to reduce to the case of split groups
and the case where $\Z_q/p^a(\!(u)\!)$
admits a pivot element.
Then apply Corollary \ref{cor:finite}.
\end{proof}

\begin{cor}
\label{cor:representable}
The morphism of the Emerton-Gee stacks
$\cX_{\lsup LG}\to\cX_{\GL_N}$
is relatively representable
by algebraic stacks of finite presentation over $\spf\Z_p$
for any embedding $\lsup LG\to\GL_N$.
\end{cor}

Here, $\cX_{\lsup LG}$
is the moduli stack of cyclotomic \'etale $(\varphi, \Gamma)$-modules
$(T, \varphi_T, \gamma_T)$
with $\lsup LG$-structure,
equipped with an identification
$T\times^{\lsup LG}\Gal(L/F)\cong
T_{\std}$,
over $\spf \Z_p$;
here $T_{\std}$
is the \'etale $(\varphi, \Gamma)$-module
with $\Gal(L/F)$-structure
that corresponds to
the quotient map $\Gal_F\to\Gal(L/F)$.

\begin{proof}
Note that $\cX_{\GL_N}\to\clR_{\GL_N}$
is relatively representable by
algebraic stacks of finite presentation over $\spf \Z_p$
(c.f. \cite[Proposition 3.4.10]{EG23}).
\end{proof}

As another application,
we establish the Shapiro's lemma.

\begin{thm} (Nonabelian Shapiro)
\label{thm:Sha}
Let $F/E$ be a finite extension.
There is an isomorphism
$\cX_{F, \lsup LG}\cong\cX_{E, \lsup L\Res_{F/E}G}$.
\end{thm}

\begin{proof}
Recall that $L$ is a splitting field of $G$.
By \cite[2.1.5]{St10},
there is a diagram
\[
\xymatrix{
H\rtimes\Gal(L/F):=(\Res_{F/E}G)^{\vee}
\rtimes \Gal(L/F)
\ar[d]\ar[r]
&
\lsup LG 
\\
\lsup L\Res_{F/E}G
}
\]

The restriction morphism
\[
\cX_{E, \lsup L\Res_{F/E}G}
\to
\cX_{F, \lsup L\Res_{F/E}G}
\]
clearly
factors through
$\cX_{F, H\rtimes\Gal(L/F)}$.
Let
$\cX_{F, H\rtimes\Gal(L/F)}
\to \cX_{F, \lsup LG}$
be the change-of-group morphism
$T\mapsto T\times^{H\rtimes\Gal(L/F)}
\lsup LG$.
So, we get a well-defined map
$f:\cX_{E, \lsup L\Res_{F/E}G}
\to \cX_{F, \lsup LG}$.

By Corollary \ref{cor:representable},
$f$ is relatively representable
by algebraic spaces of finite type over $\Z/p^a$.
By \cite[Proposition 8]{St10},
$f$ is an equivalence over artinian test objects.
The theorem now follows from \cite[Lemma 7.2.5]{LLHLM}.
\end{proof}

We don't know of a more direct proof of Shapiro's lemma;
as without the relative representability,
we cannot rule out the pathologies like
$\prod_x \spec \cX_x \to X$
where each $\cX_x$ is a local ring of a scheme $X$.

\section{The parabolic stacks $\cX_{\lsup LP}$}

In this section we assume $G$ is split
(that is, $\bG=\lsup LG$)
to simplify notation.
All results in this section can
be upgraded to general reductive groups for free.

\subsection{Parabolics}
Let $\bP=U\rtimes\bM\subset \bG$ be a parabolic subgroup
with Levi subgroup $\bM$ and unipotent radical $U$.
Let $\bB\subset \bP$ be a Borel and
$\bT\subset \bB$ be a maximal torus.

Denote by $1=U_0\subset U_1\subset \dots\subset U_n=U$
the upper central series for $U$.
For each root $\alpha\in \Phi(\bB, \bT)$,
we have a root group $\G_a\cong U_\alpha\subset \bG$.
Denote by $\Phi(U_{j}/U_i, \bT)\subset \Phi(\bB, \bT)$
the set of roots $\alpha$
such that $U_\alpha \subset U_{j}$
and $U_\alpha\not\subset U_i$,
$j>i$.

The map
(depending on an ordering of the roots
which we fix once for all)
\begin{equation}
\label{eq:splitting}
\prod_{\alpha\in \Phi(U_j/U_i, \bT)}U_\alpha
\xrightarrow{\cong}
U_j/U_i
\end{equation}
is an isomorphism of schemes.
In particular,
we have a scheme-theoretic splitting
\[
U/U_{i} \cong U/U_{i+1} \times U_{i+1}/U_i
\]
for each $i$
(which does not respect group scheme structures).

\subsection{Herr complexes}
Set 
\[
V:=U_{i+1}/U_i\xrightarrow[\cong]{\iota} \prod_{\alpha\in \Phi(U_{i+1}/U_i, \bT)}U_\alpha;
\]
since $V$ is commutative,
$\iota$ does not depend on the choice
of an ordering and is a vector space scheme
isomorphism.
Let $T^{\univ}$ be the universal
\'etale $(\varphi, \Gamma)$-module
over $\cX_{\bM}$.
Then
\[
V^\univ:=T^\univ\times^{\bM}V
\]
is an \'etale $(\varphi, \Gamma)$-module
of rank $r=\dim V$
over $\cX_{r}=\cX_{\GL_r}$.
Here $\bM\to\GL(V)$
is the adjoint action.
The Herr complex
$C^\bullet(V^\univ)$
is a perfect complex supported in degrees $[0,2]$
over $\cX_r$
(c.f. \cite[Theorem 5.1.22]{EG23})
and pulls back to a perfect complex
over $\cX_{\bM}$
and $\cX_{\bP/U_{i+1}}$.
For $\vec h=(h^1,h^2)\in \Z\times \Z$,
denote by
$\cX_{\bM}^{\vec h}\subset \cX_{\bM,\red}$
the induced reduced substack where $H^i(C^\bullet(V^\univ))$
is locally free of rank $h^i$, $i=1,2$.
Note that each 
$\cX_{\bM}^{\vec h}\subset \cX_{\bM,\red}$
is (relatively representable)
by a locally closed immersion,
by the upper-semicontinuity theorem.
Set $\cX_{\bP/U_i}^{\vec h}:=
\cX_{\bP/U_i}\underset{\cX_{\bM}}{\times}\cX_{\bM}^{\vec h}$.

\subsection{Parabolic stacks of tautological torsor type}
Denote by $\cY_{\bP/U_i}\subset \cX_{\bP/U_i}$
the fully faithful substack
such that for each $\Z/p^a$-algebra $A$,
$\cY_{\bP/U_i}(A)$ consists of \'etale
$(\varphi, \Gamma)$-modules
$(T, \varphi_T, \gamma_T)$ with $\bP/U_i$-structure whose underlying
$\bP/U_i$-torsor is isomorphic
to
$T\times^{\bP/U_i}\bM \times^{\bM} \bP/U_i$
as a $\bP/U_i$-torsor over $\A_{F,A}$.
Here $\A_{F,A}$ is the standard period ring
for cyclotomic \'etale $(\varphi, \Gamma)$-modules
(c.f. \cite{EG23} or Example \ref{ex:cyc}).
Set \[
(\bar T, \varphi_{\bar T}, \gamma_{\bar T})
:=(T, \varphi_T, \gamma_T)\times^{\bar P/U_i}\bM.
\]
Note that $\Aut_{\bP/U_i}(T)$
is a form of the algebraic group
$\bP/U_i$
and write $U'\subset \Aut_{\bP/U_i}(T)$
for the unipotent radical.
Since the formation of both
the upper central series
and of absolute root groups for
$U$ commutes with base change,
they descend to $U'$,
and we have root group decomposition
\[
U'=\prod_{\sigma\in \Phi(U, \bT)}U_\alpha'
\]
where $U_\alpha'$
is a form of $\G_a$,
and the scheme-theoretic splitting
\[
\iota_i: U'/U_{i+1}'\times U_{i+1}'/U_i'
\cong U'/U_i'
\]
also descends
as we fixed an ordering of roots once for all.

There are factorizations
\[
\xymatrix{
\varphi^*T \ar[rr]^{\varphi_{\bar T}\times^{\bM}\bP/U_i} \ar@/_3ex/[rrr]_{\varphi_T} &&T \ar[r]^{R_{X_T}} &T
}
\]
\[
\xymatrix{
\gamma^*T \ar[rr]^{\gamma_{\bar T}\times^{\bM}\bP/U_i} \ar@/_3ex/[rrr]_{\gamma_T} &&T \ar[r]^{R_{Y_T}} &T
}
\]
where $R_{(-)}$
is the right-translation by an element $(-)$ morphism
(which is the data of the underlying $\bP/U_i$-torsor
of $T$),
and $X_T, Y_T\in U'_i(\A_{F,A})$
are matrices satisfying
\[
Y_T
\Ad_{\gamma_{\bar T}}(\gamma(X_T))
=
X_T
\Ad_{\varphi_{\bar T}}
(\varphi(Y_T)).
\]
Set
$\varphi_A:=
\Ad_{\varphi_{\bar T}}\circ \varphi$
and
$\gamma_A:=
\Ad_{\gamma_{\bar T}}\circ \gamma$,
and the equation becomes
\[
Y_T\gamma_A(X_T)
=X_T\varphi_A(Y_T).
\]


\begin{lem}
The morphism $\cY_{\bP/U_i}\to \cX_{\bM}$
is relatively representable by the groupoid
whose objects over $A$ are
pairs of matrices
$(X, Y)\in U'/U'_i(\A_{F,A})$
satisfying
$Y \gamma_A(X)
=
X \varphi_A(Y)$
and whose morphisms are arrows
$(X,Y)\to (\varphi_A(Z)XZ^{-1}, \gamma_A(Z)YZ^{-1})$
for $Z\in U'/U'_i(\A_{F,A})$.
Here $(\bar T, \varphi_{\bar T}, \gamma_{\bar T})$
is the object corresponding to
$\spec A\to \cX_{\bM}$.
\label{lem:Y-to-X}
\end{lem}

\begin{proof}
The only content of the lemma
is that $\gamma_T$
is continuous if and only if $\gamma_T \times^{\bP/U_i}\bM$
is continuous.

Note that
$\bP/U_i \to \bM \times \GL(\Lie \bP/U_i)$
is an embedding.
The choice of a Levi $\bM\subset\bP$
defines
a splitting of $\Lie \bP/U_i\to \Lie \bG$,
and thus an embedding
$\GL(\Lie \bP/U_i)\to\GL(\Lie \bG)$.
So it suffices to show
$\gamma_T\times^{\bP/U_i}\Lie \bG$
is continuous,
which is equivalent to the topological nilpotency
of $(1-\gamma_T\times^{\bP/U_i}\Lie \bG)$.
We already know
$1-\gamma_T\times^{\bP/U_i}\bM \times ^{\bM}\Lie \bG$
is topologically nilpotent,
and
the difference
$\gamma_T\times^{\bP/U_i}\bM \times ^{\bM}\Lie \bG
-
\gamma_T\times^{\bP/U_i}\Lie \bG
$ is nilpotent.
\end{proof}

\subsection{Vanishing locus of higher cup products}
Next, we analyze fibers of
$\cY_{\bP/U_i}\to \cY_{\bP/U_{i+1}}$.
Again, $(\bar T, \varphi_{\bar T}, \gamma_{\bar T})$
denotes the object corresponding to
$\spec A\to \cX_{\bM}$.
Set \[
V':=U_{i+1}'/U_i'.\]
By Lemma \ref{lem:Y-to-X},
an object 
of $\cY_{\bP/U_{i+1}}$ over $\bar T$
is a pair of matrices
$X, Y\in U'/U_{i+1}'(\A_{F,A})$
and
an object 
of $\cY_{\bP/U_{i}}$ over $(X, Y)$
is a pair of matrices
$\wt X, \wt Y\in U'/U_{i}'(\A_{F,A})$
where
$\wt X=(X, \delta_X)$
and
$\wt Y=(Y, \delta_Y)$
for $\delta_X, \delta_Y\in V'(\A_{F, A})$
under the scheme-theoretic splitting
$U'/U_{i}'=U'/U_{i+1}'\times V'$.
Define
\[
o_{(X,Y)}:=
(Y,0)\gamma_A(X,0)
\varphi_A(Y,0)
(X,0)^{-1}
\in V'(\A_{F,A}).
\]
Then, we have
\[
(Y,\delta_Y)\gamma_A(X,\delta_X)
\varphi_A(Y,\delta_Y)
(X,\delta_X)^{-1}
=o_{(X,Y)}
+
(1-\varphi_A)\delta_Y
+(\gamma_A-1)\delta_X
\in V'(\A_{F,A}).
\]
The Herr complex
$C^\bullet(V'(\A_{F,A}))
=C^\bullet(\bar T\times^{\bM}V)$
is given concretely by
\[
V'(\A_{F,A})
\xrightarrow{d^0_\Herr:=(\varphi_A-1, \gamma_A-1)}
V'(\A_{F,A})
\oplus 
V'(\A_{F,A})
\xrightarrow{d^1_\Herr:=(\gamma_A-1)\oplus(1_A-\varphi)}
V'(\A_{F,A}).
\]
So,
\begin{align*}
(\delta_X,\delta_Y)
&\in C^1(\bar T\times^{\bM}V)
\\
o_{(X,Y)}
+
d^1_\Herr(\delta_X, \delta_Y)
&\in
C^2(\bar T\times^{\bM}V),
\end{align*}
and
the equation for
$(\wt X, \wt Y)$
to be an object of
$\cY_{\bP/U_i}$
is precisely
\[
o_{(X,Y)}+ d^1_\Herr(\delta_X, \delta_Y)=0.
\]
The cohomology sheaf
$H^2(C^\bullet(\bar T\times^{\bM}V))$
is a finite projective coherent sheaf
of rank $h^2$
over $\spec A$,
and thus
the short exact sequence of quasi-coherent sheaves
\[
0\to
d^1(C^1(\bar T\times^{\bM}V))
\to
C^2(\bar T\times^{\bM}V)
\to
H^2(C^\bullet(\bar T\times^{\bM}V))
\to 0
\]
admits a splitting
\[
C^2(\bar T\times^{\bM}V)
=
d^1(C^1(\bar T\times^{\bM}V))
\oplus
H^2(C^\bullet(\bar T\times^{\bM}V))
\]
Thus, we write
\[
o_{(X,Y)}
=d^1(B_X, B_Y)
+h_A
\]
where
\begin{align*}
(B_X, B_Y) &\in C^1(\bar T\times^{\bM}V)\\
h_A &\in H^2(C^\bullet(\bar T\times^{\bM}V)).
\end{align*}
The zero section of $h_A$
defines a closed subscheme of $\spec A^h\subset \spec A$;
we remark that the formation
of $\spec A^h$
does not depend on the choice of splitting,
although the element $h_A$
itself does depend on choices.

The cohomology sheaf
$H^i(C^\bullet(\bar T\times^{\bM}V))$
is a finite projective coherent sheaf
of rank $h^i$, $i=0,1$,
over $\spec A$,
and by replacing $\spec A$
by its Zariski cover,
we can assume both sheaves are finite free $A$-modules.
Let
\[
\{(z_1, w_1), \dots, (z_k, w_k)\}
\subset C^1(\bar T\times^{\bM}V))
\]
be elements whose image
in 
$H^1(C^\bullet(\bar T\times^{\bM}V))$
forms a basis of the finite free $A$-module.
We construct a vector bundle
over $\spec A^h$ by setting
\[
\cV^1_A:=
\underline{\spec}\Sym (H^1(C^\bullet(\bar T\times^{\bM}V))^\vee)
\to \spec A^h,
\]
together with a morphism
\[
\cV^1_A
\to
\cY_{\bP/U_i},
(z_i, w_i)
\mapsto
((X, z_i-B_X),(Y, w_i-B_Y)).
\]
We also define
$\cV^0_A:=
\underline{\spec}\Sym (H^0(C^\bullet(\bar T\times^{\bM}V))^\vee)$
over $\spec A^h$.

\begin{lem}
For any morphism
$\spec A\to \cY^{\vec h}_{\bP/U_{i+1}}$,
there is a Cartesian diagram
\[
\xymatrix{
[\cV^1_A/\cV^0_A] \ar[d]\ar[r]
&
\cY^{\vec h}_{\bP/U_i}\ar[d]
\\
\spec A\ar[r]
&
\cY^{\vec h}_{\bP/U_{i+1}}
}
\]
after possibly replacing $\spec A$
by a finitely presented Zariski cover.
In particular,
the
finite-type-ness of
$\cY_{\bP/U_{i+1}}^{\vec h}$
over $\Z/p^a$
implies
the finite-type-ness of
$\cY_{\bP/U_i}^{\vec h}$
over $\Z/p^a$.
\label{lem:parabolic}
\end{lem}

\begin{proof}
It is tautological
when the cohomology sheaves
are finite free over $A$.
\end{proof}

\subsection{}
We start a new subsection to clear all temporary assumptions
and notations.

\begin{cor}
There exists a finite collection
of locally closed substacks
$\{\cX_{\bM}^{\mathbf h}\}$
that exhaust $\cX_{\bM,\red}$
such that
each
$\cY_{\bP}\underset{\cX_{\bM}}{\times}\cX_{\bM}^{\mathbf h}$
is of finite presentation over $\Z/p^a$.

Moreover,
if we assume $\cX_{\bP/U_k}$ is limit-preserving
for all $k$, then we can arrange it so that
each
$
\cY_{\bP}\underset{\cX_{\bM}}{\times}\cX_{\bM}^{\mathbf h}
\to
\cX_{\bP}\underset{\cX_{\bM}}{\times}\cX_{\bM}^{\mathbf h}$
is an isomorphism.
\label{cor:parabolic}
\end{cor}

\begin{proof}
It follows from Lemma \ref{lem:parabolic}
by induction on $i$.

For the ``moreover'' part,
it suffices
to show
$\cY_{\bP/U_i}^{\vec h}
=\cX_{\bP/U_i}^{\vec h}$,
assuming
$\cY_{\bP/U_{i+1}}^{\vec h}
=\cX_{\bP/U_{i+1}}^{\vec h}$.
By Fontaine's theory
and 
the triviality
of Galois cohomology
$H^1(\bFp(\!(T)\!), U)=1$
(\cite[Proposition 6, Section III.2.1]{Se02}),
$\cY_{\bP/U_i}^{\vec h}
\underset{\cX_{\bP/U_i}^{\vec h}}{\times}\spec A\to\spec A$
is equivalence on $W(\bar\F_p)/p^a$-points
for all finite type algebras $A$
with $\spec A\to \cX_{\bP/U_i}^{\vec h}$.
A monomorphism of reduced finite type
schemes that is an equivalence on artinian points
must be an isomorphism
(c.f. \cite[Lemma 7.2.5]{LLHLM}).
Finally, since $\cX_{\bP/U_i}$
is limit-preserving,
we have
$\cY_{\bP/U_i}^{\vec h}
\cong{\cX_{\bP/U_i}^{\vec h}}$.
\end{proof}

\begin{thm}
\label{thm:rep-parabolic}
Let $G$ be an arbitray reductive group over $F$
that splits over $L$.
Let $\lsup LP=\bP\rtimes \Gal(L/F)\subset\lsup LG$
be a parabolic subgroup
with Levi factor $\lsup LM$.
Then
$\cX_{F, \lsup LP}\cong \varinjlim_m
\cX_{F, \lsup LP}^{m}$
where each
$\cX_{F, \lsup LP}^{m}$
is an algebraic stack of finite presentation
over $\spec\Z_p$
with transition maps being nilpotent thickenings.
\end{thm}

\begin{proof}
By the main theorem of \cite{Min25},
$\cX_{F, \lsup LP}$
is representable by limit-preserving
formal algebraic stacks
over $\spf\Z_p$.
By Corollary 
\ref{cor:parabolic}
and Corollary \ref{cor:representable},
$\cX_{F, \lsup LP}$ is quasi-compact.
Since $\cX_{F, \lsup LP}$ has affine diagonal,
we conclude that it is an Ind-algebraic formal algebraic stack
(\cite[Tag 0AJE]{Stacks}).
\end{proof}

\begin{remark}
We hope to conclude that
$\cX_{F,\lsup LP}\to\cX_{F, \lsup LM}$
is relatively representable by algebraic stacks
of finite presentation,
but we cannot rule out the pathology
of the following kind:
Let $\{a_1, a_2, \cdots\}$ be
a non-repeating sequence that exhausts all elements of $\bFp$.
Let $\cY$ be $\spec \bFp[t, \varepsilon]/(\varepsilon^2)$,
and let
$\cX$ be $$\dirlim{m}~\spec\bFp[t, \varepsilon]/(\varepsilon^2, (t-a_1)^m(t-a_2)^m\cdots(t-a_m)^m\varepsilon).$$
The obvious embedding
$f:\cX\hookrightarrow \cY$ is not an isomorphism
since $\cX$ is not a scheme.
Note that $\cX_{\red} = \cY_{\red}=\spec \bFp[t]$;
$f$ is fully faithful since it is an inductive limit
of closed immersions;
and $f(A)$ is clearly essentially surjective for Artinian local $\bFp$-algebras $A$.
\end{remark}

\begin{appendices}

\mychapter{Non-flat descent technique, after Emerton-Gee}
\mysection{Elementary facts about pushforward and pullback}~

Let $R$ be a Noetherian ring complete with respect
to an ideal $I\subset R$.

Consider the following
diagram
$$
\xymatrix{
    \spf \wh{R(\!(u)\!)} \ar[rd]^{\fc} & & \\
    & \spec R(\!(u)\!)\ar[r]^j & \spec R[\![u]\!] \\
    \spec \wh{R(\!(u)\!)}\ar[ru]^c & & \\
    \spec R/I^n(\!(u)\!) \ar[u]^{i_n} &
    \spec R/I^n(\!(u)\!) \ar[uu]^{i_n} &
    \spec R/I^n[\![u]\!] \ar[uu]^{i_n} &
}
$$
where $\wh{R(\!(u)\!)}$ is the $I$-adic completion of $R(\!(u)\!)$ and
the morphisms $c$, $j$, and $i_n$
are induced from the corresponding ring maps.

We collect the following simple facts.
Note that $c^*$, $j^*$, $i_n^*$ are simply base change;
$c_*$, $j_*$ and $(i_n)_*$ are forgetful functors
that preserve the underlying abelian groups;
and $\fc^*$ is by definition the $I$-adic completion functor.

\mysubsection{Fact}
(1) We have the following identities.
\begin{itemize}
\item[(a)] $j^*j_*=\id: \QCoh(R(\!(u)\!))\to \QCoh(R(\!(u)\!))$;
\item[(b)] $i_n^*(i_n)_*=\id: \QCoh(R/I^n(\!(u)\!))\to \QCoh(R/I^n(\!(u)\!))$;
\item[(c)] $i_n^*(i_n)_*=\id: \QCoh(R/I^n[\![u]\!])\to \QCoh(R/I^n[\![u]\!])$;
\item[(d)] $i^*i_*=\id: \QCoh(\wh{R(\!(u)\!)})\to \QCoh(\wh{R(\!(u)\!)})$;
\item[(e)] $\fc^*c_*=\id: \Coh(\wh{R(\!(u)\!)})\to \Coh(\wh{R(\!(u)\!)})$;
\item[(f)] $c^* = \fc^*: \Coh(R(\!(u)\!))\to \Coh(\wh{R(\!(u)\!)})$.
\end{itemize}

(2) If $M\in \Coh(R(\!(u)\!))$ and $N\in \QCoh(\wh{R(\!(u)\!)})$ is flat,
then $M\subset c_*N$ implies $c^*M\subset N$ (the map exists by adjunction).

\begin{proof}
(1-a,b,c,d) Clear.

(1-e) By \cite[Tag 00MA]{Stacks}, a finite $\wh{R(\!(u)\!)}$-module
is automatically $I$-adically complete.

(1-f) It is part (3) of \cite[Tag 00MA]{Stacks}.

(2)
By \cite[Tag 0315]{Stacks}, $\fc^* M \subset \fc^*c_*N$.
By part (1-e), $\fc^*c_*N=N$.
By part (1-f), $\fc^*M=c^*M$.
\end{proof}


%

\mysection{Mittag-Leffler systems of modules}~

See \cite[Tag 0594]{Stacks} for the basic definitions of Mittag-Leffler systems.

\begin{dfn}
Let $(A_i, \varphi_{ij})$, $(B_i, \varphi_{i,j})$
be Mittag-Leffler directed inverse systems of $R$-modules.
We say
\begin{itemize}
\item $(A_i, \varphi_{ij}) \subset (B_i, \varphi_{ij})$
if $A_i\subset B_i$ for each $i$;
\item $(A_i, \varphi_{ij}) \models (B_i, \varphi_{ij})$
if $(A_i, \varphi_{ij}) \subset (B_i, \varphi_{ij})$
and $\varphi_{ji}(A_j)=\varphi_{ji}(B_j)$ for all $i$
and all $j\gg i$.
\end{itemize}
\end{dfn}

\begin{lem}
Let $(A_i, \varphi_{ij})$, $(B_i, \varphi_{ij})$ be
Mittag-Leffler directed inverse systems of $R$-modules,
and let $(D_i, \varphi_{ij})$ be a system that contains both
$(A_i, \varphi_{ij})$, $(B_i, \varphi_{ij})$.
Then $(A_i \cap B_i, \varphi_{ij})$ is Mittag-Leffler.

Moreover if $(C_i, \varphi_{ij})\subset (D_i, \varphi_{ij})$
and $(B_i, \varphi_{ij})\models (C_i, \varphi_{ij})$,
then
 $(A_i\cap B_i, \varphi_{ij})\models (A_i\cap C_i, \varphi_{ij})$.
\end{lem}

\begin{proof}
Clear.
\end{proof}

\begin{lem}
For each $t\in R$,
if $(A_i, \varphi_{ij})$ is a
Mittag-Leffler directed inverse system of $R$-modules
then so is $(A_i\otimes_RR[1/t], \varphi_{ij})$.

Moreover if $(A_i, \varphi_{ij})\models (B_i, \varphi_{ij})$,
then
 $(A_i\otimes_RR[1/t], \varphi_{ij})\models (B_i\otimes_RR[1/t], \varphi_{ij})$.
\end{lem}

\begin{proof}
It suffices to show that
if $f:M\to N$ is a homomorphism of $R$-modules,
then $\Img(f)\otimes_R R[1/t] = \Img(f\otimes_R R[1/t])$.
Note that $M\to f(M)$ is surjective and $f(M)\subset N$ is injective.
After localization, $M[1/t]\to f(M)[1/t]$ is surjective
and $f(M)[1/t]\to N[1/t]$ is injective since $R[1/t]$ is $R$-flat.
So $f(M)[1/t]$ is the image of $f[1/t]$.
\end{proof}

\begin{lem}
Fix an integer $c>0$.
Let $(A_i, \varphi_{ij})$, $(B_i, \varphi_{ij})$ be
Mittag-Leffler directed inverse systems of $R$-modules
where indexes $i\in \Z$.
Assume there exists surjective homomorphisms
$f_i: A_i\to B_i$, $g_i: B_i\to A_{i-c}$
for all $i$ such that
\begin{itemize}
\item Both $f_i$, $g_i$ are compatible with transition maps $\varphi_{ij}$;
\item $g_i\circ f_i = \varphi_{i,i-c}$.
\end{itemize}
Then $\invlim{i} A_i = \invlim{i} B_i$.
\end{lem}

\begin{proof}
Clear since both $\lim f_i \circ \lim g_i$
and $\lim g_i \circ \lim f_i$ are identities.
\end{proof}

\mysection{Non-flat descent of vector bundles}~

\begin{dfn}
\label{def:fun-ext-vect}
Given a tuple $(U, X, Y, j, \pi)$
where $j:U\to X$ is a scheme-theoretically dominant
flat morphism of affine schemes
and $\pi:Y\to X$ is a scheme-theoretically dominant
proper morphism of schemes.
Denote by $\Vect_{U, X, Y}$
the category of triples
$(V, K, \theta)$ where $V$ is a locally free coherent sheaf over $U$,
$K$ is a locally free coherent sheaf over $Y$,
and $\theta$ is an isomorphism $\theta: \pi^*V \cong j^*K$.
For ease of notation, we will drop $\theta$
and write $(V, K, \theta)=(V,K)$ and $\pi^*V = j^*K$.
\end{dfn}

Write $Y_U:= U\times_X Y$.
We have a diagram
$$
\xymatrix{
    Y_U \ar[d]^\pi \ar[r]^j & Y \ar[d]^\pi\\
    U \ar[r]^j & X
}
$$
where all morphisms are scheme-theoretically dominant.

\begin{lem}
\label{lem:fun-ext-1}
We have $\pi_*j_*\pi^* = \pi_*\pi^*j_*$
as functors from coherent sheaves over $U$ to quasi-coherent sheaves
over $X$.
\end{lem}

\begin{proof}
Since the maps involved are either
proper or affine,
pushforward preserve quasi-coherence;
since $X$ is affine,
to show
$\pi_*j_*\pi^*F = \pi_*\pi^*j_*F$
it suffices to check the global sections.
Note that by the pullback-pushforward adjunction,
we have a canonical map $\pi^*j_*\to j_*\pi^*$.
Since $j^*\pi_*j_*\pi^*F$
is in set-theoretic bijection with
$\pi_*j_*\pi^*F$,
it suffices to show
$j^*\pi_*j_*\pi^*F = j^*\pi_*\pi^*j_*F$
which is clear since $j^*j_*=1$ and $j^*\pi_*=\pi_*j^*$
by flat base change.
\end{proof}

By the previous lemma,
for each $(V, K)\in \Vect_{U,X,Y}$,
    $\pi_*\pi^*j_*V = j_*j^*\pi_*K$.
Since $\pi$ and $j$ are both scheme-theoretically
dominant, the units of the pullback-pushforward adjunctions
$1\to \pi_*\pi^*$ and $1\to j_*j^*$
are inclusion of quasi-coherent sheaves.

\begin{dfn}
Define a functor
$(j, \pi)_*: \Vect_{U, X, Y} \to \QCoh_X$
by $(V, K)\mapsto j_*V \cap \pi_*K$ ($\subset \pi_*\pi^*j_*V$).
\end{dfn}

\begin{lem}
$(j, \pi)_*: \Vect_{U, X, Y}\to \QCoh_X$
is a lax monoidal functor (see Appendix \ref{tan-cat} for definitions).
\end{lem}

\begin{proof}
Let $(V_1, K_1)$, $(V_2, K_2)\in \Vect_{U, X, Y}$.
Since both $j_*$ and $\pi_*$ are lax monoidal functors,
we have maps
\begin{align*}
(j_*V_1) \otimes (j_*V_2) &\to j_*(V_1\otimes V_1) \subset \pi_*\pi^*j_*(V_1\otimes V_2)\\
(\pi_*K_1) \otimes (\pi_*K_2) &\to \pi_*(K_1\otimes K_1) \subset \pi_*\pi^*j_*(V_1\otimes V_2)
\end{align*}
whose restriction to $(j, \pi)_*(V_1, K_1)\otimes (j,\pi)_*(V_2, K_2)$
coincide.
The restriction defines a lax monoidal structure on $(j, \pi)_*$.
\end{proof}

\begin{lem}
\label{lem:fun-ext-2}
If $j$ is an open immersion, $j^*(j, \pi)_*(V, K)=V$.
\end{lem}

\begin{proof}
We have $j_*j^*(j_* V\cap \pi_*K) = j_*V \cap j_*j^*\pi_*K
=j_*V\cap j_*\pi_*j^*K = j_*V \cap \pi_*\pi^*j_*V=j_*V$.
Note that $j^*j_*=\id$.
\end{proof}

Write $Y[\![u]\!]$ for $Y\times_{\spec R} \spec R[\![u]\!]$.

\begin{prop}
\label{prop:decomp}
Let $R = \invlim{n} R/I^n$ be a Noetherian complete ring.
Let $\pi:X\to \spec R$ be
a scheme-theoretically dominant, proper morphism.
There exists a lax monoidal functor
$$
\xi:\Vect_{\spec \wh{R(\!(u)\!)}, \spec{R[\![u]\!]}, X[\![u]\!]}
\to \Coh_{\spec R(\!(u)\!)}
$$
which becomes the forgetful functor
after composing with the base change functor
$\Coh_{\spec R(\!(u)\!)}\to \Coh_{\spf \wh{R(\!(u)\!)}}$.
\end{prop}

\begin{proof}
Write $R_n$ for $R/I^n$, and
write $X_n$ for $X\times_{\spec R}\spec R_n$.
The arrows in the following diagram
$$
\xymatrix{
    X_n(\!(u)\!) \ar[d]^\pi \ar[r]^j & X_n[\![u]\!] \ar[d]^\pi\\
    \spec R_n(\!(u)\!) \ar[r]^j & \spec R_n[\![u]\!]
}
$$
are all scheme-theoretically dominant.
Write $i_n: \spec R_n \to \spec R$
for the embedding.
Write $j$ for $\spec R(\!(u)\!)\to \spec R[\![u]\!]$.
For $(V, K)\in \Vect_{\spec \wh{R(\!(u)\!)}, \spec{R[\![u]\!]}, \mathfrak{X}[\![u]\!]}$,
set $$\kappa_n(V, K) := (j, \pi)_*i_n^*(V,K)$$ and
$\kappa := \invlim{n} \kappa_n$.
Define $\xi(V, K):=j^* \kappa(V,K).$

By the theorem on formal functions,
$\invlim{n} \pi_*i_n^*K = \pi_*K$ is coherent and thus the submodule
$\kappa(V,K)$ is also coherent.
So $\xi(V, K)\in \Coh(R(\!(u)\!))$.

It is clear that $\xi$ is a lax monoidal functor.
So it remains to show $\xi$ composed with
the $I$-adic completion functor is the forgetful functor;
in other words,
$$
\text{($\dagger$)}\hspace{5mm}\invlim{n} i_n^*\xi (V, K)
=\invlim{n} j^* i_n^* \kappa(V, K)= V.
$$

By the Artin-Rees lemma (\cite[Tag 00IN]{Stacks}),
there exists a constant number $c>0$ such that
$$
(I^n \pi_*K)\cap \kappa(V, K) = I^{n-c}(I^c \pi_*K\cap \kappa(V, K))
$$
for all $n\ge c$.
So we have
\begin{align*}
I^n \kappa(V,K) & \subset (I^n \pi_*K) \cap \kappa(V,K) \\
(I^{n+c} \pi_*K) \cap \kappa(V,K) & \subset I^n \kappa(V,K)
\end{align*}
for all $n$.
Note that
$$
\frac{\kappa(V,K)}{(I^n \pi_*K) \cap \kappa(V,K)}
=\Img(\kappa(V,K)\to i_n^*\pi_*K) =:A_n.
$$
The following maps
\begin{align*}
A_n &\to i_n^* \kappa(V,K) \\
i_n^* \kappa(V,K) &\to A_{n+c}
\end{align*}
are surjective for all $n$.
It follows that
$$
\invlim{n} j^* i_n^*\kappa(V,K) = \invlim{n} j^* A_n.
$$
By \cite[Tag 02OB]{Stacks}, there exists a positive integer
$c'$ such that
if we write $K_n:=\Img(\pi_*K \to \pi_*i_n^*K)$,
then $K_n \to i^*_{n-c'}\pi_*K$
is surjective for all $n\ge c'$.
Write $B_n:=\Img(\kappa(V,K) \to \pi_*i_n^*K)\subset K_n$.
We have surjective maps $A_n \to B_{n-c'}$ for all $n \ge c'$.
Note that
$(B_n) \models (\kappa_n(V,K))$. So
$$
\invlim{n} j^* \kappa_n(V,K) = \invlim{n} j^* B_n.
$$
By Lemma \ref{lem:fun-ext-2},
$$
\invlim{n} j^* \kappa_n(V,K) = \invlim{n} i_n^*V = V.
$$
Note that both $(A_n)$ and $(B_n)$ are Mittag-Leffler inverse systems with
surjective transition maps,
and that $\invlim{n} A_n =\invlim{n} B_n= \kappa(V,K)$;
it is possible to choose an infinite sequence of increasing integers
$m_1 < m_2 < \dots$
such that $A_{m_i}$ is a quotient of $B_{m_{i+1}-c'}$ for all $i$.
Now it is clear that
$$
\invlim{n} j^* A_n=
\invlim{i} j^* A_{m_i}=
\invlim{i} j^* B_{m_i-c'}=
\invlim{i} j^* B_{n}
=V.
$$
So we are done.
\end{proof}

\mychapter{Infinite-dimensional descent theory, after Drinfeld}

\mysection{Descending $\cG$-torsors}
\label{ssec:descent-torsor}

Let $R$ be a ring which is not necessarily noetherian.

\begin{lem}
\label{tan-lem1}
(1) Let $M_1\to M_2\to M_3$ be a sequence of $R(\!(u)\!)$-modules.
Let $R\to S$ be an fpqc ring map.
Assume $0\to M_1\otimes_{R(\!(u)\!)}S(\!(u)\!)\to M_2\otimes_{R(\!(u)\!)}S(\!(u)\!)\to
M_3\otimes_{R(\!(u)\!)}S(\!(u)\!)\to 0$ is a short exact sequence
of finitely generated projective $S(\!(u)\!)$-modules,
then $0\to M_1\to M_2\to M_3\to 0$ is a short exact sequence
of finitely generated projective $R(\!(u)\!)$-modules.

(2) The same is true for $R[\![u]\!]$.
\end{lem}

\begin{proof}
(1)
By \cite[Theorem 5.1.18]{EG21}, each of $M_i$ is finitely generated projective.
By Lemma \ref{lem-eff2}, it suffices to show
that if $f:\spec S\to \spec R$ is surjective,
then each closed point of $\spec R(\!(u)\!)$ is contained in the image of $f$.
Let $\fp$ be a maximal ideal of $R(\!(u)\!)$.
Then $\fp_0 := \{\text{Leading coefficient of $f$}|f\in \fp\}$
is a proper ideal of $R$.
Since $\spec S\to \spec R$ is surjective, $\fp_0S\ne S$.
Hence $\fp S(\!(u)\!)\ne S(\!(u)\!)$, and
$\spec S(\!(u)\!)\times_{\spec R(\!(u)\!)}\spec R(\!(u)\!)/\fp$
is non-empty.

(2) The proof is completely similar to that of (1).
\end{proof}

Let $\cG$ be a smooth affine group scheme of finite type over $\cO$
whose connected components all admit an integral point.
For an $\cO$-scheme $X$,
let $\gtor_X$ be the category of $\cG$-torsors over $X$.

Let $\frepg(\cO)$ be the category of 
linear representations of $\cG$ on finitely generated projective $\cO$-modules.

We recall the following theorem of \cite[Theorem 2.1.1]{Lev15}.

\begin{lem}
\label{tan-thm0}
Let $X$ be an algebraic space over $\spec\cO$.
Then there is an equivalence of categories
$F:\gtor_X \xrightarrow{\sim} [\frepg(\cO), \vect_X]^\otimes$,
which commutes with arbitrary base change on $X$.
\end{lem}

\begin{proof}
\cite[2.1.1]{Lev15} or \cite[2.5.2]{Lev13} proves this when
$\cG$ is a group with connected geometric fibres and when $X$ is a scheme.
However, the connectivity requirement is only used to guarantee that
$\cO_{\cG}$ is $\cO$-projective.

The algebraic space case follows formally from the scheme case by descent.
\end{proof}

We fix some notations.
Let $P$ be a $\cG$-torsor over $X$.
Let $r: \cG \to \GL(V)$ be a representation of $\cG$.
Set $P\times^{\cG}V := F(P)(r)$.
Also denote by $P\times^{\cG}\GL(V)$ the 
principal bundle associated to
$P\times^{\cG}V$.

\begin{lem}
\label{tan-lem-2}
Let $R\to S$ be an fpqc ring map.

(1)
The category of $\cG$-torsors over $\spec R(\!(u)\!)$
is equivalent to the category of $\cG$-torsors over $\spec S(\!(u)\!)$
with descent datum.

(2)
The category of $\cG$-torsors over $\spec R[\![u]\!]$
is equivalent to the category of $\cG$-torsors over $\spec S[\![u]\!]$
with descent datum.
\end{lem}

\begin{proof}
Follows formally from Lemma \ref{tan-lem1} and Theorem \ref{tan-thm0}.
\end{proof}

\mysection{The $\cF$-twisted affine Grassmannian}
\label{taGR}

This main result of this section is an elaboration of the remark under Propositon 3.8
of \cite{Dri06}.
The terminology ``$\cF$-twisted affine Grassmannian''
also comes from \cite{Dri06}.

Let $R$ be an arbitrary ring.
Let $\cG$ be a smooth affine $\cO$-group scheme of finite type.

We fix some notations.
Write $D_R$ for $\spec R[\![u]\!]$ and $D_R^\star$ for $\spec R(\!(u)\!)$.
Let $X$ be an $R$-scheme and let $R\to S$ be a ring homomorphism.
Write $X\otimes_RS$ for $X\times_{\spec R}\spec S$.

Let $P$ be a $\cG$-torsor over $D_R^\star$.
Define the groupoid $\GR_P$ over $\spec R$ as follows.
Let $R\to S$ be a ring homomorphism.
Set $\GR_P(S)$ to be the groupoid of pairs
$(T,\gamma)$ where $T$ is a $\cG$-torsor over $D_S$
and $\gamma$ is an isomorphism 
$T|_{D_S^\star}\xrightarrow{\sim} P\times_{D_R^\star}D_S^\star$.
Note that when $P = \cG_R := \cG \otimes_{\cO} R$ is the trivial $\cG$-torsor,
$\GR_{\cG_R}$ is the usual affine Grassmannian.

As a consequence of Lemma \ref{tan-lem-2},
the groupoid $\GR_{P}$ is stacky.
Moreover, because of the framing, $\GR_{P}$ is a setoid.

\begin{lem}
When $\cG=\GL_N$,
the groupoid $\GR_P$ 
is representable by an ind-proper algebraic space.
\end{lem}

\begin{proof}
It is \cite[Proposition 3.8]{Dri06}.
\end{proof}

\mysubsection{}
Let $R\to S$ be a ring homomorphism.
By unwinding the definitions,
there is an isomorphism of groupoids
$\GR_{P\times_{D_R^\star}D_S^\star}\xrightarrow{\sim} \GR_P\times_{\spec R}\spec S$.

Let $i: \cG\to \GL_N$ be a closed immersion, then there is
an obviously defined pushforward morphism
$i_*:\GR_P \to \GR_{i_* P}$, where
$i_* P:= P\times^{\cG}\GL_N$.

\begin{lem}
\label{gr-lem1}
Assume $\cG$ be a connected reductive group scheme over $\cO$.
Let $R\to S$ be a ring homomorphism.
Then the map
$|\GR_P(S)| \to |\GR_{i_*P}(S)|$
is injective.
\end{lem}

\begin{proof}
Suppose $|\GR_P(S)|\ne \emptyset$.
Let $(T,\gamma)\in \GR_P(S)$.
Since $\cG$ is smooth, there is an \'etale cover $\spec S'\to \spec S$
such that $T\times_{D_S}D_{S'}$ is a trivial $\cG$-torsor over $D_{S'}$.
Now that $\GR_P\times_{\spec R}\spec S'\to \GR_{i_*P}\times_{\spec R}\spec S'$
is relatively representable by a closed immersion
(see, for example, \cite[Corollary 3.3.10]{Lev13}).
By descent (\cite[Tag 04SK]{Stacks} and \cite[Tag 0420]{Stacks}),
$\GR_P\times_{\spec R}\spec S\to \GR_{i_*P}\times_{\spec R}\spec S$
is also relatively representable by a closed immersion.
Hence
$|\GR_P(S)| \to |\GR_{i_*P}(S)|$
is injective.
\end{proof}

\begin{lem}
\label{gr-lem2}
The morphism of groupoids $\GR_P\to \GR_{i_*P}$
is relatively representable by a closed immersion.
\end{lem}

\begin{proof}
Fix $(M,\gamma_M)\in \GR_{i_*P}(R)$.
Write $X$ for $\GR_{P}\times_{\GR_{i_*P}}\spec R$.
An object of $X(S)$ is a pair $((T, \gamma), \alpha)$
where $(T,\gamma)\in \GR_P(S)$ and $\alpha$
is an isomorphism $i_*T \xrightarrow{\sim} M\otimes_{R[\![u]\!]}S[\![u]\!]$
which is compatible with $\gamma$ and $\gamma_M$.

Define the groupoid $Y$ over $R[\![u]\!]$-algebras as follows:
set $Y(B)$ to be the groupoid of pairs $(T, \alpha)$
where $T$ is a $\cG$-torsor over $\spec B$ and 
$\alpha$ is an isomorphism $i_*T \xrightarrow{\sim}M\otimes_{R[\![u]\!]}B$.
Let $R\to S$ be a ring homomorphism.
There is a map
$X(S)\to Y(S[\![u]\!])$ which sends
$((T,\gamma),\alpha)\mapsto (T,\alpha)$.
By Lemma \ref{gr-lem1},
the map $|X(S)|\to |Y(S[\![u]\!])|$ is injective.
It is standard that $Y$ is represented by an affine scheme
(namely, by $M\times^{\GL_N}(\GL_N/\cG$)).
The image of $X(S)$ are objects $(T,\alpha)$ of $Y(S[\![u]\!])$
whose restriction on $D_S^\star$
is $(P\times_{D_R^\star}D_S^\star, (P\to i_*P)\times_{D_R^\star}D_S^\star)$.
By \cite[Lemma 3.3.9]{Lev13},
there exists a closed subscheme $\spec R/I \subset \spec R$ such that
$|X(S)|$ is non-empty if and only
$R\to S$ factors through $\spec R/I$.
Hence $X \cong X\times_{\spec R}\spec R/I$.
By replacing $R$ by $R/I$,
we can assume $X(R) \ne \emptyset$.
Now that there is an fppf ring map $R\to S$ such that
$P\times_{D_R^\star}D_S^\star$ is a trivial $\cG$-torsor;
and $X\times_{\spec R}\spec S \cong \GR_{\cG}\times_{\GR_{\GL_N}}\spec S$,
which is well-known to be a closed subscheme of $\spec S$.
This lemma is now proved by descent.
\end{proof}

\begin{lem}
\label{gr-lem3}
Let $X$ be a closed subscheme of the ind-scheme $\GR_P$.
There is a Nisnevich cover $\spec S\to \spec R$ such that
$X\times_{\spec R}\spec S$ is a projective scheme over $\spec S$.
\end{lem}

\begin{proof}
By Lemma \ref{gr-lem2}, it is reduced to the $\GL_N$-case (\cite[Proposition 3.8]{Dri06}).
\end{proof}

\mychapter{Tannakian categories}
\label{tan-cat}

An \textit{exact category} is an additive category where a class
of short exact sequences is specified.
An \textit{exact functor} is an additive functor
which takes a short exact sequence
to a short exact sequence.

A \textit{monoidal category} is a tuple $(\cC,-\otimes-, \bI)$
where $\cC$ is a category,
$-\otimes-:\cC\times \cC\to \cC$ is a bifunctor called the tensor product,
and $\bI$ is an object of $\cC$ called the unit object.
For a pair of objects $(X, Y)$, the \textit{internal Hom} 
$\iHom(X,Y)$ (if it exists) is defined to be
the object representing the functor
$T \mapsto \Hom(T\otimes X, Y)$.
For an object $X$, the \textit{dual} $X^\vee$ (if it exists)
is defined to be the object representing $\iHom(X, \bI)$.

Let $\cC$, $\calD$ be monoidal categories.
A \textit{lax monoidal functor} from $\cC$ to $\calD$
is a pair $(F, \alpha)$ where $F$ is a functor from
$\cC$ to $\calD$ sending $\bI$ to $\bI$,
and $\alpha$ is a natural transformation
from the bifunctor $F(-)\otimes F(-)$ to 
the bifunctor $F(-\otimes -)$
satisfying some coherence conditions.
$F$ is said to be a \textit{strict monoidal functor}
if $\alpha$ is a natural isomorphism.
If a strict monoidal functor has a right adjoint functor,
then the right adjoint functor is canonically a lax monoidal functor.

Let $\cC$ be an exact monoidal category.
An object $X$ is said to be an \textit{invertible object}
if the functor $-\otimes X$ is an exact equivalence.

Let $\calD$, $\cE$ be exact monoidal categories.
Denote by $[\calD, \cE]^{\otimes}$
the category of faithful, exact, strict monoidal functors
$\calD\to\cE$.
A morphsim $\alpha:F\to G$ is a natural transformation satisfying
$\alpha_{A\otimes B}=\alpha_A\otimes\alpha_B$ for $A, B\in\calD$
and $\alpha_{\bI}=\id$.

A \textit{rigid category} is a monoidal category $C$ such that
(1) internal Hom always exists; (2) the morphism
$\iHom(X_1, Y_1)\otimes\iHom(X_1,Y_2) \to \iHom(X_1\otimes X_2,Y_1\otimes Y_2)$ is an isomorphism for objects $X_1$, $Y_1$, $X_2$ and $Y_2$;
(3) the morphism $X \to (X^{\vee})^{\vee}$ is an isomorphism.
A strict monoidal functor between rigid categories automatically
preserves inner Hom and duality \cite[Proposition 1.9]{DM82}.

A monoidal category is said to be a \textit{symmetric monoidal category}
if it is equipped with a natural isomorphism
$s_{A,B}: A\otimes B \to B\otimes A$ satisfying the obvious
coherence conditions.

\mysection{}

Let $Z$ be a scheme.
The category $\Coh_Z$ of coherent sheaves on $Z$
is an abelian symmetric monoidal category.
The category $\Vect_Z$ of finitely generated projective $\cO_Z$-modules
is an exact symmetric monoidal category.
Let $f: Z\to Y$ be a proper morphism.
The pullback functor $f^*:\Coh_Y \to \Coh_Z$
is a strict monoidal functor,
and hence its right adjoint $f_*: \Coh_Z \to \Coh_Y$
is a lax monoidal functor.
\footnote{
    When $f: \spec B \to \spec A$ is a morphism of affine schemes,
    the lax monoidal structure of $f_*$ is given by the map
    $M\otimes_AN \to M\otimes_B N$.
}

\mysection{Setting and some monoidal categories}
\label{et-phi-mod}
Let $R$ be a ring.
Let $Z$ be an $R$-scheme.
Let $U \subset Z$ be an open dominant subscheme.
Let $\varphi: Z\to Z$ be a morphism which map $U$ to $U$.
\begin{itemize}
\item 
Define $\Coh_{Z,U,\varphi}$ to be the category of
coherent $\cO_Z$-modules $\cF$, together with a homomorphism
$\phi:\varphi^*\cF|_U \to \cF|_U$
(which we call a $\varphi$-structure).
\item
Define $\Coh_{Z,U,\varphi}^\et$ to be the full subcategory
of $\Coh_{Z,U,\varphi}$ consisting of objects whose $\varphi$-structure
is an isomorphism.
\item
Define $\Vect_{Z,U,\varphi}$ (resp. $\Vect_{Z,U,\varphi}^\et$) to be the full subcategory
of $\Coh_{Z,U,\varphi}$ (resp. $\Coh_{Z,U,\varphi}^\et$) consisting of finitely generated projective
$\cO_Z$-modules.
Note that $\Vect_{Z, U, \varphi}^{\et, \ef}$
is also a full subcategory of $\Vect_{Z, \varphi}$.
\end{itemize}
We also call objects of $\Coh_{Z,U,\varphi}^{\et}$ {modules
with \'etale $\varphi$-structure}.

\begin{lem}
\label{cat-lem2}
The category $\Vect_{Z,U,\varphi}^\et$ is 
an exact, rigid, symmetric monoidal category.
\end{lem}

\begin{proof}
Since the inverse image functor $\varphi^*$ is a strict monoidal functor,
the tensor product of two \'etale $\phi$-structure
is an \'etale $\phi$-structure.
The unit object in $\Vect_{Z,U,\varphi}^\et$
is the structure sheaf $\cO_Z$ together with the identification
$\varphi^* \cO_Z  = \cO_Z\otimes_{\varphi^{-1}\cO_Z}\cO_Z = \cO_Z$.

Moreover, since the inverse image functor $\varphi^*$
preserves the sheaf Hom, internal Homs in $\Vect_{Z,U,\varphi}^\et$
are representable.
\end{proof}

\begin{lem}
\label{lem:lG-tan}
For each $\Z_p$-scheme $X$,
the category of $\lG$-torsors on $X$
is equivalent to the category of fiber functors
from $\frep_{\lsup LG}$ to $\Vect_X$.
\end{lem}

\begin{proof}
It is \cite[Theorem 2.5.2]{Lev13}.
Note that the connectedness assumption is only used
in \cite[Proposition C.1.8]{Lev13}.
Since $\lG$ as a scheme is a disjoint union of $\wh G$,
\cite[Proposition C.1.8]{Lev13} holds for $\lG$.
\end{proof}

\mychapter{Exactness of Frobenius decompletion}~

We will need the following descent theorem
which strengthens \cite[Theorem 5.5.20]{EG21}.

\begin{lem}
\label{lem-eff1}
Let $R$ be a complete Noetherian local $\F_p$-algebra with
maximal ideal $\fm$,
and let $\varphi: R(\!(u)\!)\to R(\!(u)\!)$
be the $R$-linear map that sends $u$ to $u^q$,
where $q=p^f$.
\end{lem}

Let $M_1\to M_2\to M_3$ be a sequence of finite free
\'etale $\varphi$-modules over $R(\!(u)\!)$.
If $0\to \widehat{M_1}\to\widehat{M_2}\to \widehat{M_3}\to 0$
is a short exact sequence
(where $(-)^{\wedge}$ denotes the $\fm$-adic completion),
then $0\to M_1\to M_2\to M_3\to 0$ is also an exact sequence.

\begin{proof}
Let $S$ be a ring which contains an element $u$.
Denote by $S\langle - \rangle$ the $u$-adic completion of $S[-]$.

By the Cohen structure theorem,
$R = K[[T_1, \dots, T_n]]/I$ for some field extension $K$ of 
$\F_p$
and some ideal $I$.
For each integer $m \ge 0$, we define
$$
R_{m} := R[u]\langle \frac{T_1}{u^{q^m}},\dots, \frac{T_n}{u^{q^m}}\rangle[1/u]
$$
and
$$
R_{-m} := R[u]\langle \frac{T_1^{q^m}}{u},\dots, \frac{T_n^{q^m}}{u}\rangle[1/u].
$$
It is helpful to visualize $\spec R_m$ ($-\infty < m < \infty$) as disks of various radii.
\begin{center}
\begin{tikzpicture}
\coordinate (O) at (0,0);

\coordinate (O) at (0,0);
\draw[fill=red!30] (O) circle (3.4);
\draw[fill=green!40] (O) circle (2.6);
\draw[fill=yellow!70] (O) circle (1.8);
\draw[fill=white!45] (O) circle (1);
\node at (0, 0) {\dots};
\draw[decoration={text along path,reverse path,text align={align=center},text={Spec{$R_{m+1}$}}},decorate] (1.3,0) arc (0:180:1.3);
\draw[decoration={text along path,reverse path,text align={align=center},text={Spec{$R_m$}}},decorate] (2.1,0) arc (0:180:2.1);
\draw[decoration={text along path,reverse path,text align={align=center},text={Spec{$R_{m-1}$}}},decorate] (2.9,0) arc (0:180:2.9);
\end{tikzpicture}
\end{center}

We advise the readers to read \cite[5.5.20]{EG21} for
the following facts:
\begin{itemize}
\item[(i)]
There are flat homomorphisms of $R(\!(u)\!)$-algebras
$R_m \hookrightarrow R_{m+1}$ determined by $T_i\mapsto T_i$ 
and $R(\!(u)\!)$-isomorphisms
$$
R(\!(u)\!)\otimes_{\varphi,R(\!(u)\!)}R_m\xrightarrow{\sim} R_{m+1}
$$
for $m = \dots,-1,0,1,\dots$;
\item[(ii)]
Write $\varphi_m: R_m\to R_{m+1}$ for the composition
$R_m\xrightarrow{x\mapsto 1\otimes x}R(\!(u)\!)\otimes_{\varphi,R(\!(u)\!)}R_m\xrightarrow{\sim} R_{m+1}$,
which is a faithfully flat map;
\item[(iii)]
The ring $R_{\infty}:=\dirlim{m\ge 0}R_m$ is faithfully flat over 
$\widehat{R(\!(u)\!)}$, the $\fm$-adic completion of $R(\!(u)\!)$;
\item[(iv)]
Each maximal ideal $\fp$ of $R(\!(u)\!)$ comes from
a maximal ideal of $R_{m}$ for some $m\le 0$.
More precisely, if $R(\!(u)\!)\to L$ is a surjection
onto a field $L$ (which is necessarily a finite extension of $K(\!(u)\!)$),
then it factors through a surjection $R_m\to L$ for some $m\le 0$.
\end{itemize}

We will use the following constructions.

\begin{itemize}
\item[(I)] For each $R((U))$-module $M$, write $j_m^*M$ for $R_m\otimes_{R(\!(u)\!)}M$, which can be interpreted geometrically as restriction to the ``disk'' $\spec R_m$.
\item[(II)] For each $R_m$-module $M$, write $\varphi_m^*$ for
$R_{m+1}\otimes_{\varphi_m, R_m}M$,
which can be interpreted geometrically as the Frobenius amplification of $M$.
\end{itemize}
The following fact is crucial:
$$\varphi_m^*j_m^*M = j_{m+1}^*\varphi^*M.$$

Since $R_{\infty}$ is faithfully flat over $\widehat{R(\!(u)\!)}$,
the sequence
$$0\to M_1\otimes_{R(\!(u)\!)} R_{\infty}
\to M_2\otimes_{R(\!(u)\!)} R_{\infty}
\to M_3\otimes_{R(\!(u)\!)} R_{\infty}\to 0$$
is exact.
it is possible to choose a basis $\{b_1,b_2,\dots,b_s;c_1,c_2,\dots,c_t\}$
of $M_2\otimes_{R(\!(u)\!)}R_{\infty}$ such that 
the image of $M_1\otimes_{R(\!(u)\!)} R_{\infty}$ is generated by $\{b_1,\dots,b_s\}$
and each element of $M_3\otimes_{R(\!(u)\!)}R_{\infty}$
can be lifted to an element of the submodule generated by $\{c_1,\dots,c_t\}$.
Since $R_{\infty} = \dirlim{m}R_m$,
there exists some $N\ge 0$ such that
$\{b_1,b_2,\dots,b_s;c_1,c_2,\dots,c_t\}\subset M_2\otimes_{R(\!(u)\!)}R_N$
and as a consequence,
the sequence
$$0\to j_N^*M_1
\to j_N^*M_2
\to j_N^*M_3\to 0$$
is exact.

Since $\varphi: R(\!(u)\!)\to R(\!(u)\!)$ is faithfully flat and thus
the base change $\varphi_m: R_m\to R(\!(u)\!)\otimes_{\varphi,R(\!(u)\!)}R_m
\cong R_{m+1}$ is also faithfully flat.
We have a commutative diagram
$$
\xymatrix{
    0 \ar[r] & j_m^* M_1 \ar[r] &
    j_m^* M_2 \ar[r] &
    j_m^* M_3 \ar[r] & 0\\
    0 \ar[r] & j_m^* \varphi^*M_1 \ar[r]\ar[u]^{\phi_{M_1}}\ar@{=}[d] &
    j_m^* \varphi^*M_2 \ar[r]\ar@{=}[d]\ar[u]^{\phi_{M_2}} &
    j_m^* \varphi^*M_3 \ar[r]\ar@{=}[d]\ar[u]^{\phi_{M_3}} & 0 \\
    0 \ar[r] & \varphi_{m-1}^*j_{m-1}^* M_1 \ar[r] &
    \varphi_{m-1}^*j_{m-1}^*M_2 \ar[r] &
    \varphi_{m-1}^*j_{m-1}^*M_3 \ar[r] & 0
}
$$
where all vertical maps are isomorphisms.
Hence the short exactness of 
$$0\to j_m^*M_{1}\to j_m^*M_{2}\to j_m^*M_{3}\to 0$$
implies the short exactness of
$$0\to j_{m-1}^*M_{1}\to j_{m-1}^*M_{2}\to j_{m-1}^*M_{3}\to 0.$$

To show $0\to M_1\to M_2\to M_3$ is exact,
it suffices to show that for each maximal ideal $\fp$
of $R(\!(u)\!)$, the localisation at $\mathfrak {p}$ of this sequence is exact.
We have already mentioned that there exists some integer $m\le 0$
and some maximal ideal $\fq$ of $R_{m}$ which
lies above $\fp$.
Since $R(\!(u)\!)\to R_{m}$ is flat,
$R(\!(u)\!)_{\fp}\to (R_m)_{\fq}$ is faithfully flat.
The proof is now complete.
\end{proof}

The following elementary lemma allows us to generalise the above lemma
to a general $\Z_p$-algebra $A$.

\begin{lem}
\label{lem-eff2}
Let $A$ be a ring.
Let $f: \spec B\to \spec A$ be a morphism such that
all the closed points of $\spec A$ are contained in the image of $f$.
Let $M_1\to M_2\to M_3$ be a sequence of finite projective $A$-modules
such that $$0\to M_1\otimes_AB\to M_2\otimes_AB\to M_3\otimes_AB\to 0$$ is
a short exact sequence.
Then $$0\to M_1\to M_2\to M_3\to 0$$ is also a short exact sequence.
\end{lem}

\begin{proof}
Since the property of being a short exact sequence is local,
it suffices to consider the case where both $A$ and $B$ are local rings,
and $A\to B$ is a local homomorphism.
We can replace $B$ by its residue field.
So assume $A$ is a local ring
with maximal ideal $I$, and
$B=A/I$.

Fix a spliting
$\bar\psi:M_1/IM_1 \oplus M_3/IM_3\xrightarrow{\sim} M_2/IM_2$
of the short exact sequence.
Let $\psi:M_1\oplus M_3 \to M_2$ be an arbitrary lift of $\bar\psi$
which is compatible with the given lift of the short exact sequence;
since $M_2\to M_3$ is surjective by Nakayama's lemma,
such a lift $\psi$ exists.
By Nakayama's lemma again,
$\psi$ is a surjective homomorphism of finite free $A$-modules
of the same rank,
and is thus an isomorphism (by a determinant argument).
\end{proof}

\begin{lem}
\label{lem-eff3}
Let $R$ be a complete Noetherian local
$\Z_p$-algebra killed by $p^a$ for some integer $a\ge 1$,
with maximal ideal $\fm$,
and let $\varphi:R(\!(u)\!)\to R(\!(u)\!)$
be a $R$-linear deformation of the $q$-power Frobenius
such that $\varphi(u)\subset R[\![u]\!]$.
Let $M_1\to M_2\to M_3$ be a sequence of finitely generated
\'etale $\varphi$-modules with $R$-coefficients.

If the $\fm$-adic completion of $0\to M_1\to M_2\to M_3\to 0$
is short exact, then $0\to M_1\to M_2\to M_3\to 0$
is a short exact sequence of projective $\Lambda_R(\!(u)\!)$-modules.
\end{lem}

\begin{proof}
By \cite[Theorem 5.5.20]{EG21}, $M_i$ is a projective module, $i=1,2,3$;
note that the running assumption of \cite{EG21}
is $\varphi(R[\![u]\!])\subset R[\![u]\!]$.
In particular, each $M_i\otimes \F_p$
is projective over $R/p(\!(u)\!)$. 
By the last two paragraphs
of the proof of \cite[Theorem 5.5.20]{EG21}
the exactness in the mod $p$ case implies
exactness in the general case.
By \cite[Lemma 5.2.14]{EG21},
one can find finite projective \'etale $\varphi$-modules 
with $R$-coefficients
$N_1$ and $N_3$ such that both $M_1\oplus N_1$ and $M_3\oplus N_3$
are finite free.
By Lemma \ref{lem-eff1} and Lemma \ref{lem-eff2},
$0\to M_1\oplus N_1\to M_2\oplus N_1\oplus N_3\to M_3\oplus N_3\to 0$
is a short exact sequence.
\end{proof}

\begin{lem}
\label{lem-eff4}
Let $\cC$ be an exact monoidal category,
and let $R$ be a complete
Noetherian local $\Z/p^a$-algebra
with defining ideal $I$.
Let $F$ be a lax monoidal functor from $\cC$
to $\Vect_{\spec R(\!(u)\!)}$
such that the $I$-adic completion $\wh F\in
[\cC, \Vect_{\spec \wh{R(\!(u)\!)}}]^\otimes$ of $F$
a faithful, exact, strict monoidal functor.
If $F$ is equipped with an \'etale $\varphi$-structure,
then
$F$ is a faithful, exact, strict monoidal functor.
\end{lem}

\begin{proof}
Faithfulness of $F$ follows from the faithfulness of $\widehat F$.
Exactness of $F$ is Lemma \ref{lem-eff3}.
Lemma \ref{lem-eff3} also
implies that $F$ is strict monoidal.
\end{proof}

\end{appendices}

\printbibliography

\end{document}